\documentclass[10pt]{article}

\usepackage[latin1]{inputenc}
\usepackage{import}
\usepackage[pagebackref=true, colorlinks=true, linkcolor=blue, citecolor=blue, linkcolor=blue, anchorcolor=red, urlcolor=blue]{hyperref}
\usepackage{url, amsmath,enumerate,fancyhdr,amssymb, amsthm, 
pstricks, epsf, lscape, ifpdf, graphicx, float}
\usepackage{algorithm}
\usepackage{enumitem}
\usepackage{color}
\usepackage{algpseudocode} 
\algnewcommand\algorithmicinput{\textbf{Input:}}
\algnewcommand\INPUT{\item[\algorithmicinput]}
\algnewcommand\algorithmicinitialization{\textbf{Initialization:}}
\algnewcommand\INITIALIZATION{\item[\algorithmicinitialization]}

\newcommand{\pha}{\phantom} 
\newcommand{\ti}{\times}

\newcommand{\Black}[1]{\textcolor{black}{#1}} 
\newcommand{\mrbmath}[1]{\mbox{{\Red{{\boldmath ${#1}$}}}}}
\newcommand{\mblbmath}[1]{\mbox{{\Black{{\boldmath ${#1}$}}}}}

\newcommand{\mbbmath}[1]{ \mbox{\Blue{{\boldmath ${#1}$}}}}

\newcommand{\mgbmath}[1]{ \mbox{\Green{{\boldmath ${#1}$}}}}
\newcommand{\mblackbmath}[1]{ \mbox{\Black{{\boldmath ${#1}$}}}}

	\newcommand\bovermat[2]{%
		\makebox[0pt][l]{$\smash{\overbrace{\phantom{%
						\begin{matrix}#2\end{matrix}}}^{\text{#1}}}$}#2}
	
	\usepackage{mathtools}
	\makeatletter
	\renewcommand*\env@matrix[1][*\c@MaxMatrixCols c]{%
		\hskip -\arraycolsep
		\let\@ifnextchar\new@ifnextchar
		\array{#1}}
	\makeatother
	
	\usepackage{tikz,array,calc, xcolor}
	\usetikzlibrary{decorations.pathreplacing}

\newcommand{\Red}[1]{\textcolor{red}{#1}}
\newcommand{\Blue}[1]{\textcolor{blue}{#1}}
\newcommand{\Green}[1]{\textcolor{green}{#1}}

\usepackage[margin=2.7cm]{geometry}

\renewcommand{\arraystretch}{1.2}

\mathchardef\mhyphen="2D

\newtheorem{Definition}{Definition}

\newtheorem{Example}{Example}
\newtheorem{Proposition}{Proposition}
\newtheorem{Lemma}{Lemma}
\newtheorem{Theorem}{Theorem}
\newtheorem{Corollary}{Corollary}
\newtheorem{Remark}{Remark}
\newtheorem{Assumption}{Assumption}

\newcommand{\ba}{\begin{array}}  
	\newcommand{\ena}{\end{array}}

\newcommand{\A}{\mathcal A}

\newcommand{\const}{\operatorname{constant}}
\newcommand{\trace}{\operatorname{trace}}

\newcommand{\bpx}{\begin{pmatrix}}
\newcommand{\epx}{\end{pmatrix}}
\newcommand{\bbx}{\begin{bmatrix}}
\newcommand{\ebx}{\end{bmatrix}}

\newcommand{\bdef}{\begin{Definition}} 
\newcommand{\commentout}[1]{}
\newcommand{\co}[1]{}

\newcommand{\coab}[1]{}

\newcommand{\nin}{\noindent}

\newcommand{\pf}[1]{\vspace{.35cm} \nin {\bf Proof {#1} }}

\newcommand{\sym}[1]{{\cal S}^{#1}}
\newcommand{\psd}[1]{{\cal S}_+^{#1}}
\newcommand{\rad}[1]{\mathbb{R}^{#1}}

\newcommand{\symn}{\sym{n}}
\newcommand{\psdn}{\psd{n}}

\newcommand{\I}{ {\cal I} }

\definecolor{mygreen}{RGB}{0,120,0}

\newcommand{\beq}{\begin{equation}}
\newcommand{\eeq}{\end{equation}}
\newcommand{\beqa}{\begin{eqnarray}}
\newcommand{\eeqa}{\end{eqnarray}}
\newcommand{\bac}{\begin{array}{ccccccccccc}}
\newcommand{\eac}{\end{array}}
\newcommand{\bprop}{\begin{Proposition}}
\newcommand{\eprop}{\end{Proposition}}

\newcommand{\beqast}{\begin{eqnarray*}}
\newcommand{\eeqast}{\end{eqnarray*}}
\newcommand{\benum}{\begin{enumerate}}
\newcommand{\eenum}{\end{enumerate}}
\newcommand{\bit}{\begin{itemize}}
\newcommand{\eit}{\end{itemize}}
\newcommand{\bth}{\begin{Theorem}}
\newcommand{\enth}{\end{Theorem}}
\newcommand{\ble}{\begin{Lemma}}
\newcommand{\ele}{\end{Lemma}}
\newcommand{\bex}{\begin{Example}}
\newcommand{\eex}{\end{Example}}
\newcommand{\bcor}{\begin{Corollary}}
\newcommand{\ecor}{\end{Corollary}}
\newcommand{\brem}{\begin{Remark}}
\newcommand{\erem}{\end{Remark}}
\newcommand{\bass}{\begin{Assumption}}
\newcommand{\eass}{\end{Assumption}}

\renewcommand{\arraystretch}{1.2}

\newcommand{\bsmx}{\begin{small} \begin{pmatrix}}
\newcommand{\esmx}{\end{pmatrix} \end{small}}

\title{\Large How do exponential size solutions arise in semidefinite programming?}  
\author{G\'{a}bor  Pataki, Aleksandr Touzov \thanks{Department of Statistics and Operations Research, University of North Carolina at Chapel Hill} \hspace{1cm} 
}

\begin{document}

\maketitle

\begin{abstract}
	A striking pathology  of semidefinite programs (SDPs) is   illustrated by a classical example of Khachiyan: feasible solutions in SDPs may need exponential space  even to write down. 
		Such exponential size solutions are  the main obstacle to solve a long standing, fundamental  open problem: can we decide feasibility of SDPs in polynomial time? 
	
		The consensus seems that SDPs with large size solutions are rare. 
		However, here we prove that they are actually quite common: a linear change of variables transforms every strictly feasible SDP into a Khachiyan type SDP, in which the leading variables 
	are large. As to ``how large", that depends on the singularity degree of a dual problem.
	Further, we present some SDPs coming from sum-of-squares proofs,  in which 
	large solutions  appear naturally, without any change of variables. 
	We also partially answer the question: how do we represent such large solutions in polynomial space?

\end{abstract}



{\em OR/MS subject classification:} Primary: convexity; secondary: programming-nonlinear-theory

\section{Introduction} \label{sect-intro}

\tableofcontents

\paragraph{Linear programs and polynomial size solutions} 
The classical  linear programming (LP) feasibility problem asks whether a system of linear inequalities 
$$
x_1 a_1 + \dots + x_m a_m + b \geq 0 
$$
has a solution, where the $a_i$ and $b$ are column vectors with integer entries. 
 When the answer is ``yes", 
 then by a classical argument there exists  a feasible rational $x$ whose  size  is  at most $2 m^2 \log m$ times the size of 
 the matrix $[a_1, \dots, a_m, b].$ 
 When the answer is ``no",  there is a certificate of infeasibility whose size is similarly bounded. 
 
 Here and in the sequel  we define 
 size (or bit-length)  as in  \cite[Section 2.1]{schrijver1998theory}.   Precisely, the size of a rational number $p/q, \, $ where $p$ and $q$ are relatively prime integers, is $\lceil \log_2 (| p| + 1) \rceil + \lceil \log_2 (| q| + 1 ) \rceil +1.$
  The size of a  rational vector with $k$ elements is the sum of the sizes of its elements plus $k;$ and 
  the size of an $k \times \ell$ rational matrix is the sum of the sizes of its elements plus $k \cdot \ell.$ The size of 
  a rational number/vector/matrix is essentially  the number of bits needed to describe it in binary representation. 
  
\paragraph{Semidefinite programs and exponential size solutions} 
Semidefinite programs (SDPs) are a far reaching generalization of linear programs, and in recent decades
 they have attracted widespread  attention. An SDP feasibility problem can be formulated as
\begin{equation} \label{problem-sdp} \tag{\mbox{$\mathit{P}$}} 
x_1 A_1 + \dots + x_m A_m+ B \succeq 0, 
\end{equation} 
where the $A_i$ and $B$ are symmetric matrices  with integer entries. We assume that the $A_i$ are linearly independent, and as usual,  $S \succeq 0$  means that the symmetric matrix $S$ is positive semidefinite. 

 In striking contrast to a linear program, in some cases all feasible solutions of \eqref{problem-sdp} have  exponential size  
 as a function of the number of variables.  This  surprising fact is illustrated by a classical  convex feasibility problem 
 of  Khachiyan:
\begin{equation} \label{problem-khachiyan} \tag{\mbox{\em{Khachiyan}}} 
x_1 \geq x_2^2, \, x_2 \geq x_3^2, \dots, x_{m-1} \geq x_m^2, \, x_m \geq 2.
\end{equation}
We first show how fast the variables grow in \eqref{problem-khachiyan},
	so suppose $x$ is feasible in  it. 
	Then by a straightforward calculation we get
	$x_1 \geq 2^{2^{m-1}},
	$
	so $\log_2 x_1 \geq 2^{m-1}.$ Hence the size of $x_1, \, $ and of any feasible solution  is at least $2^{m-1}.$

We next show how to cast \eqref{problem-khachiyan} in the form of \eqref{problem-sdp}, so 
we  write its 
quadratic constraints as 
\begin{equation} \label{eqn-khach-2by2} 
\bpx x_i & x_{i+1} \\ x_{i+1} & 1 \epx \succeq 0 \,\, \text{for} \,\, i=1, \dots, m-1.
\end{equation}
Then we define a symmetric matrix $\A(x)$ with $2m - 1$ rows (and columns) as follows. We set the first $m-1$ two by two
principal blocks of $\A(x)$ 
equal to the matrices in \eqref{eqn-khach-2by2}, and the lower right corner of $\A(x)$ equal to $x_m - 2.$
Then $\A(x) \succeq 0$ holds if and only if $x$ satisfies \eqref{problem-khachiyan}. Finally, it is straightforward to put the problem $\A(x) \succeq 0$ into the form of \eqref{problem-sdp} by defining suitable $A_i$ and $B$ matrices. 

\begin{figure}[h]
\begin{center}
	\includegraphics[scale=0.5]{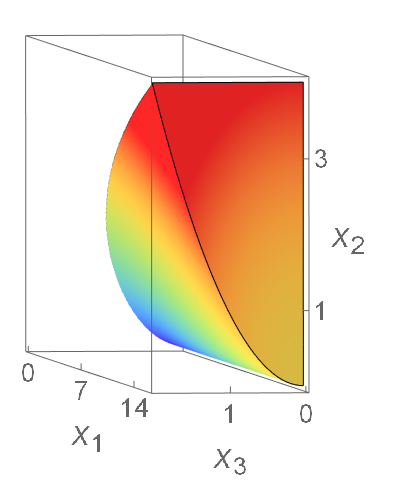}
	\includegraphics[scale=0.5]{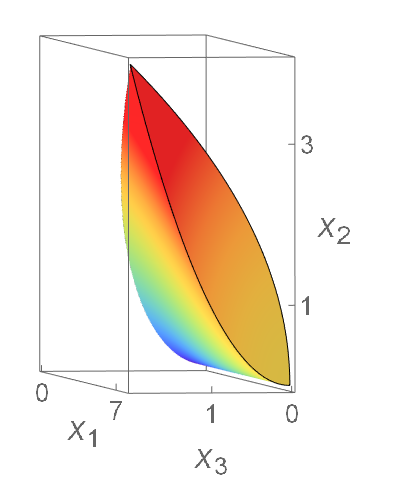}
\end{center}
\caption{Feasible sets of \eqref{problem-khachiyan} (on the left) and of the quadratic inequalities \eqref{eqn-mild-quadratic-for-pic} derived from  \eqref{problem-mild-sdp}  (on the right). } 
	\label{figure-together} 
\end{figure}

We show the feasible set of  \eqref{problem-khachiyan}  with $m=3  \, $ on the left  in 
Figure \ref{figure-together}. Our goal is to illustrate how fast $x_1$ and $x_2$ grow with respect to $x_3. \, $ 
Hence, to better visualize this growth rate (and not run out of space),  we replaced the constraint $x_3 \geq 2$ by $2 \geq x_3 \geq 0$ and made 
$x_3$ increase from right to left. 

	Observe that Khachiyan's example shows much more than exponential size solutions (in the number of variables) 
	 in an SDP may {\em exist}. 
	After all, even in a linear program with unbounded feasible set 
	solutions of any size exist!   
	However, in \eqref{problem-khachiyan} {\em all} solutions must have exponential size; and for that, the key is the hierarchy among the variables. 

To be precise, 
from now on an 
``SDP with exponential size solutions" will mean 
an SDP in which {\em all} feasible solutions have exponential size in the number of variables.

Why are we interested in SDPs with exponential size solutions? Mostly because they  are the main obstacle to solving 
the following fundamental open problem: 
	$$
	\text{Can we decide feasibility of \eqref{problem-sdp} in polynomial time? }
	$$
	Indeed,  algorithms that decide feasibility of \eqref{problem-sdp} 
	in polynomial time  must assume that a polynomial size solution exists (if there is  a solution to start with): 
	see a detailed exposition in  \cite{grotschel2012geometric}.  Algorithms that optimize a linear function over the set defined by \eqref{problem-sdp}
	in polynomial time also need similar assumptions 
	\cite{nesterov1994interior, Ren:01, de2016turing}.  
	In  contrast, the algorithm in  \cite{porkolab1997complexity} that achieves 
	the best known complexity bound  to decide feasibility of SDPs  uses  a fundamental result  from  the first order theory of reals
	\cite{renegar1992computational},   and it runs in polynomial time only in fixed dimension.
	
We know of few papers that deal directly with the complexity of SDP. However, several works study the complexity of a related problem, optimizing a polynomial subject to polynomial inequality constraints. 
On the positive side, some polynomial optimization problems are 
polynomial time  solvable when the dimension is fixed: see 
\cite{renegar1992computational, bienstock2016note, bienstock2020complexity, barvinok1993feasibility,  vavasis1990proving}.
Further,  polynomial size solutions exist in special cases \cite{vavasis1990quadratic}.  
 On the other hand, several fundamental problems in polynomial optimization are NP-hard, see for example, \cite{bienstock2020complexity, pardalos1991quadratic,  ahmadi2013np, ahmadi2019complexity}.

					Khachiyan's example inevitably  leads to   the following questions:  
			\begin{eqnarray} 
				  \label{question-exp-common} 
				& \text{Do SDPs with exponential size solutions occur frequently? } & \\
				 \label{question-exp-represent} 
				& \text{In such SDPs, can we represent  the feasible solutions   in polynomial space} ? 
					\end{eqnarray} 
			

The answer to \eqref{question-exp-common} seems to be a ``no",  since 
the only such SDP we know of is \eqref{problem-khachiyan}. 
However, to question \eqref{question-exp-represent} we have hope to get  a ``yes" answer.
After all, 
to convince ourselves 
that 
$x_1 := 2^{2^{m-1}}$ (with a suitable $x_2, \dots, x_m$)
is feasible in \eqref{problem-khachiyan}, we do not need to write down
$x_1$ explicitly:  
instead, we can just do a symbolic computation. 
Still, question \eqref{question-exp-represent}  seems to be open.

		\paragraph{Contributions}  
Perhaps surprisingly, 
	we will  
answer ``yes" to question \eqref{question-exp-common}; and we will give a partial ``yes" answer to  question \eqref{question-exp-represent}. 
One of the underlying techniques we use is facial reduction  \cite{BorWolk:81, Pataki:00B, Pataki:13, drusvyatskiy2017many, WakiMura:12}
that was originally introduced to 
induce strong duality in conic optimization problems.

				We assume that \eqref{problem-sdp} has a {\em strictly feasible} solution $x$ for which  
				$\sum_{i=1}^m x_i A_i + B$ is positive definite.
											We fix a nonnegative integer parameter $k, \, $ the  {\em singularity degree} of a dual problem.
				We will precisely define $k$  
				soon, but for now we only need to   know two facts. First,  $k \leq 1$ holds when \eqref{problem-sdp} is a linear program; and second, that $k=m$ holds  in the 
				SDP representation of \eqref{problem-khachiyan}. 
			
			An informal version of our main result follows.
			\co{ \footnote{Theorem \ref{thm-main} (and its informal version) would hold vacuously even when $k \leq 1,$ so we 
				assume $k \geq 2$ for simplicity. 
		}			}
				
			{\em 	{\bf Informal Theorem \ref{thm-main}}
		Suppose $k \geq 2.$ Then there is an invertible matrix $M$ such that the linear change of variables  
		$x \leftarrow Mx \, $ transforms 
			\eqref{problem-sdp}  into a problem 	\eqref{problem-sdpprime} with the following properties:
			
			If $x$ is strictly feasible in \eqref{problem-sdpprime}, and $x_k$ is sufficiently large, then $x_1, x_2, \dots, x_k$ 	obey a Khachiyan type hierarchy.  Precisely,  the inequalities
	\begin{equation} \label{eqn-xj-xj+1-pre}
		x_1 \geq d_2 x_2^{\alpha_2}, \, x_2 \geq d_3  x_3^{\alpha_3},  \dots, \, x_{k-1}  \geq d_k  x_k^{\alpha_k}  
	\end{equation} 
	hold, where 
	\begin{equation} \label{eqn-alphaj-bounds-pre} 
		2  \geq \alpha_{j}  \geq  1 + \frac{1}{k-j+1} \,\,\,  \text{for} \, j=2, \dots, k. 
	\end{equation}
	Here the $d_j$ and 	$\alpha_j$ are positive constants 
	that depend on the $A_i, \,$ on  $B, \, $ and the last $m-k$ variables, that we consider fixed. }
\qed

Loosely speaking, Theorem \ref{thm-main} can be interpreted as:  ``a linear transformation 
uncovers a Khachiyan type hierarchy  in all strictly feasible SDPs."

Suppose $x$ is as stated in Informal Theorem \ref{thm-main}, in particular $x_k$ is large and positive. 
	Then  $x_1$ is larger than $x_k.$

How much larger?  To find out, we combine the inequalities \eqref{eqn-xj-xj+1-pre} and \eqref{eqn-alphaj-bounds-pre}. 
First we look at the worst case, when $\alpha_j =2 $ for all $j, \, $  just like 
in  \eqref{problem-khachiyan}.  Then  $x_1$ is at least constant times $x_k^{2^{k-1}},$ so the size of $x_1$ is exponentially larger than that of $x_k.$ 
On the other hand, suppose   $\alpha_j = 1 + \frac{1}{k-j+1}$ for all $j. \, $  Then 
by an elementary calculation we get that even in this best case $x_1$ is at least constant times $x_k^{k}.  \, $ 

We next discuss why our assumptions to prove Theorem \ref{thm-main} are minimal. 
First, we must assume  that 
\eqref{problem-sdp} has a strictly feasible solution.   Indeed, 
there are SDPs without strictly feasible solutions, with large dual singularity degree, but 
no Khachiyan type hierarchy among the variables, and no large solutions. We discuss such an SDP after 	Example \ref{example-khanch-sdp}.

Next we explain why in Theorem \ref{thm-main}  we must allow 
a linear change of variables.  
	Suppose we perform a linear 
	change of variables, say $x \leftarrow Gx $ in \eqref{problem-khachiyan} ,  
	 where $G$ is a random, dense matrix. 
After this change  \eqref{problem-khachiyan} will  be quite messy, and 
will have no  variables that are obviously larger than others.. Thus,  if we performed this transformation,  then we must perform its inverse   $x \leftarrow G^{-1}x$ to get back to \eqref{problem-khachiyan}.
Finally, we argue that in Theorem \ref{thm-main} we must focus 
	on just a subset of variables and restrict the last of these variables 
	to be sufficiently large. Indeed, suppose we 
	replace the constraint $x_{m} \geq 2$ by $x_m \geq 2 + x_{m+1} \, $  in \eqref{problem-khachiyan}, where $x_{m+1}$ is a new variable. 
	After this change $x_1, \dots, x_m$ can all be zero, so there is no longer a  hierarchy among them, and none of them is  forced to be large.
	Thus, we must restrict  $x_m$ to be larger than $1 \, $ 
	 (even though $x_m > 1$ is now not implied by the constraints) to restore the hierarchy. 
	
	Besides proving Theorem \ref{thm-main}, we show that in SDPs coming from minimizing a univariate polynomial a Khachiyan type hierarchy, and
		large variables appear naturally; that is, without a change of variables, and without assuming that 
		$x_k$ is large enough. 
		The same is true of an SDP published in \cite{o2017sos} that proves nonnegativity of a linear function 
		over a set described by quadratic constraints.

We will also partially answer the representation question \eqref{question-exp-represent} as follows.
Inequalities \eqref{eqn-xj-xj+1-pre} and \eqref{eqn-alphaj-bounds-pre} imply that whenever $x$ is strictly feasible in the transformed SDP  \eqref{problem-sdpprime}, variables $x_1, \dots, x_k$ can take 
on large values. 
However,  we will see that to verify that a strictly feasible $x$ in \eqref{problem-sdpprime} exists, we will never have to compute these large values numerically. Instead, we will just do a symbolic computation to convince ourselves that suitable values of $x_1, \dots, x_k$ exist.
See the discussion after the proof of Lemma \ref{lemma-xk-arb-large}.

\paragraph{Related  work}

Linear programs can be solved in polynomial time, as it was first proved by Khachiyan \cite{khachiyan1979polynomial}; see Gr{\"o}tschel, Lov{\'a}sz, and Schrijver  \cite{grotschel2012geometric} 
 for an exposition that handles important details like the necessary accuracy. 
Other landmark polynomial time algorithms for linear programming were given   by Karmarkar 
\cite{karmarkar1984new}, Renegar \cite{renegar1988polynomial}, and Kojima et al \cite{kojima1989primal}.

On the other hand,  to decide SDP feasibility  in polynomial time, we must assume that  there is a 
polynomial size solution (provided there is a solution). We refer to  
\cite{grotschel2012geometric} for such an algorithm based on the ellipsoid method.
The algorithm of Porkolab and Khachiyan \cite{porkolab1997complexity} is the fastest known algorithm 
to decide  SDP feasibility; however, it  runs in polynomial time only for fixed $n$ and $m.$ The  algorithm of \cite{porkolab1997complexity} 
uses a foundational  result of 
Renegar \cite{renegar1992computational}, which decides in polynomial time the 
feasibility of a system of polynomial inequalities in fixed dimension. 
We further refer to Nesterov and Nemirovskii \cite{nesterov1994interior} for foundational  interior point methods to solve SDPs with an objective function.
We also refer to Renegar \cite{Ren:01} for a very clean  treatment of interior point methods for convex optimization; and to DeKlerk and Vallentin \cite{de2016turing} 
for a very precise bit complexity analysis of interior point methods  to solve  SDPs. 

The complexity of  SDP is closely related to the complexity of optimizing a polynomial subject to polynomial inequality 
constraints. To explain how, first  consider a 
system of convex quadratic inequalities 
\begin{equation} \label{problem-qcqp} 
x^\top Q_i x + b_i^\top x + c_i \leq 0 \; (i=1, \dots, m)
\end{equation}
where the $Q_i$ are fixed symmetric psd matrices, and $x \in \rad{n}$ is the vector of variables. The question whether 
we can decide feasibility of \eqref{problem-qcqp} in polynomial time is also fundamental, open,  
and, in a sense, easier than the question of deciding feasibility of \eqref{problem-sdp} in polynomial time. 
The reason is that \eqref{problem-qcqp} can be represented as an instance of \eqref{problem-sdp} by choosing suitable 
$A_i$ and $B$ matrices. 
On the other hand, 
 we can formulate semidefiniteness of a symmetric matrix variable by requiring the principal minors (which are polynomials) 
 to be nonnegative.  

Among positive results in polynomial optimization,  we already mentioned Renegar's paper \cite{renegar1992computational}.
\co{by a foundational result of Renegar \cite{renegar1992computational} 
we can decide feasibility of a system of polynomial
inequalities in fixed dimension. This result is the basis of the algorithm in \cite{porkolab1997complexity}.}
Bienstock \cite{bienstock2016note} proved that such problems can be solved in polynomial time, if the number of 
constraints is fixed, the constraints and objective are quadratic,
and at least one constraint is strictly convex. Further, Sakaue et al \cite{sakaue2016solving} designed a practical algorithm to solve 
such problems with two constraints. The work of \cite{bienstock2016note}
builds on Barvinok's fundamental result \cite{barvinok1993feasibility} that proved we can test in polynomial time
whether  a system of a fixed number of  quadratic equations is feasible. 
It also builds on 
 early work of Vavasis \cite{vavasis1990quadratic} which proved that a system with 
linear constraints  and one quadratic constraint has a solution of polynomial size. In other important early work, Vavasis and Zippel \cite{vavasis1990proving} proved  we can solve indefinite quadratic optimization problems with a ball constraint, in polynomial time.
Other related papers are  e.g., by Stern and Wolkowicz \cite{stern1995indefinite}, and Pong and Wolkowicz 
\cite{pong2014generalized}. These  show that the trust region subproblem with an indefinite objective function can be viewed as a convex problem, 
and hence solved efficiently.

On the flip side, there are many hardness results. For example, Bienstock, del Pia, and Hildebrand 
\cite{bienstock2020complexity} proved  it is NP-hard to test whether a system of quadratic inequalities 
has a polynomial size rational solution, even if we know that the system has a rational solution. 
Pardalos and Vavasis \cite{pardalos1991quadratic} proved the fundamental problem of minimizing a (nonconvex) quadratic function subject to linear constraints  is also NP-hard. The following problem is also classical, and was proven to be NP-hard only in 2013, by 
Ahmadi, Olshevsky, Parrilo, and Tsitsiklis \cite{ahmadi2013np}:  can we test convexity of a polynomial? 
It is also NP-hard to test whether a polynomial optimization problem  attains its optimal value, see Ahmadi and Zhang \cite{ahmadi2019complexity}. 

One of the tools we use is an elementary facial reduction algorithm.
These algorithms were originally designed to ensure  strong duality in conic optimization problems. 
They originated in the paper of Borwein and Wolkowicz \cite{BorWolk:81}, then simpler variants were given,  for example, by Waki and Muramatsu \cite{WakiMura:12} and in \cite{Pataki:00B, Pataki:13}. 
	For a recent comprehensive survey of facial reduction and its applications, see Drusvyatskiy and Wolkowicz \cite{drusvyatskiy2017many}. 

In other related work, O' Donnell 
\cite{o2017sos} presented an SDP that certifies  nonnegativity of a polynomial via the sum of squares (SOS) proof system, and is essentially equivalent to \eqref{problem-khachiyan}.  Previously it was thought that sum-of-squares  proofs, a popular tool in theoretical computer science,  can be found in polynomial time.
However, due to O' Donnell's work, it is now  clear that this is not obviously the case. Precisely, the complexity of finding SOS proofs is just as open as the complexity of 
deciding feasibility of SDPs. 


\paragraph{The plan of the paper}
		 In Subsection \ref{subsection-prelim} we review preliminaries.
		In Subsection \ref{subsection-state-thm1} we formally state Theorem \ref{thm-main} and illustrate it via two extreme examples. 
		In Subsection \ref{subsection-proofs} we  prove it  
		 in a sequence of lemmas. In particular, in 
		Lemma \ref{lemma-xj>=dj+1} we give a recursive formula, akin to a continued fractions formula, 
		to compute the $\alpha_j$ exponents in \eqref{eqn-xj-xj+1-pre}. As an alternative, in Subsection \ref{subsection-fourier-motzkin} we show how to compute the $\alpha_j$  using the classical  Fourier-Motzkin elimination for linear inequalities; this is an interesting contrast, since 
		SDPs are highly nonlinear. 	 In Section  \ref{section-polyopt} we cover the case of SDPs coming from polynomial optimization and also revisit the example from 
		\cite{o2017sos}.  Section 
			 \ref{section-conclusion} concludes with a discussion.

			 Our proofs are fairly elementary.  We use Proposition \ref{prop-gordan-stiemke},  a convex analysis argument 
			 about positive semidefinite matrices and linear subspaces. However, other than that, 
			 we only rely on basic linear algebra, and on manipulating quadratic polynomials.

			\subsection{Notation and preliminaries} \label{subsection-prelim} 
			
			
			\paragraph{Matrices} 
			Given a matrix $M \in \mathbb{R}^{n \times n}$ and 
			$R, S  \subseteq \{1, \dots, n\}$ we denote the submatrix of $M$ corresponding to rows in $R$ and columns in $S$ by 
			$M(R,S). $ 	We write $M(R)$ to abbreviate  $M(R,R).$ 
			\co{ When $R = \{r\}$ and $S = \{s\}$ are singletons, we simply write 
			$M(r,s)$ to denote $M(\{r\}, \{s\})$ and $M(r)$ to denote $M(r,r).$ 
			} 
			
			We let $\symn$ be the set of $n \times n$ symmetric matrices and $\psdn$ be the set of 
			$n \times n$ symmetric positive semidefinite (psd) matrices. The notation $S \succ 0$ means that the symmetric matrix $S$ is positive definite. 
		The  inner product of symmetric matrices $S$ and $T$ is defined 
			as $S \bullet T := \trace(ST).$

				\begin{Definition} \label{definition-regfr} 
			 We say that $(C_1, \dots, C_\ell)$ is a {\em regular facial reduction sequence for $\psd{n}$} 
				if each $C_i$ is in $\symn$ and of the form 
						$$
					C_1  = 
					\bordermatrix{
						& \overbrace{\qquad}^{\textstyle r_{1}} &  \overbrace{\qquad \qquad}^{\textstyle{n-r_1 }} \cr\\
						& I   &  0   \cr
						& 0   &  0  \cr}, \, \dots, C_i  = 
				\bordermatrix{
					& \overbrace{\qquad \qquad \qquad}^{\textstyle r_{1}+\ldots+r_{i-1}} & \overbrace{\qquad}^{\textstyle r_{i}} & \overbrace{\qquad \qquad \qquad\quad}^{\textstyle n-r_{1}-\ldots-r_{i}} \cr\\
					& \times  &  \times  &  \times \cr
					& \times  &  I   &  0 \cr
					& \times  &  0  &  0 \cr} 
				$$
				for  $i= 1,\dots, \ell, \,$ where the $r_i$ are nonnegative integers, and the $\times$ symbols 
				correspond to blocks with arbitrary elements.
				
					\end{Definition}
					To provide background, we next explain the parlance ``facial reduction sequence." For that, suppose $Y$ is a psd matrix, which has zero $\bullet$ product with $C_1, \dots, C_\ell.$  
					Since $C_1 \bullet Y = 0, \, $ the sum of the first $r_1$ diagonal elements of $Y$ is zero. Since these diagonal elements are nonnegative,  they must be all zero. Thus, since $Y \succeq 0, \, $  its first $r_1$ rows and columns are zero.  Hence, $C_2 \bullet Y \, $ is the sum of the diagonal elements of $Y$ in rows 
				 $r_1+1, \dots, r_1+ r_2$ and we similarly deduce that these rows (and corresponding columns) of $Y$ 
				 are all zero. 
					
					Continuing, we learn that  $Y$ is reduced to live in the set 
					$$
					F \, = \, \{ \, Y \in \psd{n} \, : \, \text{the first} \,\, r_1+ \dots + r_\ell  \, \text{rows and columns of } \, Y \, \text{are zero} \,  \},
					$$
					and we know that $F$ is a face of $\psdn$ 
				\footnote{A convex subset $F$ of  $\psdn$ is a face, if for 
			any $X, Y \in \psdn$ if the open line segment $\{ \lambda X + (1-\lambda)Y :   0 < \lambda < 1 \}$ intersects $F, \, $ then both $X$ and $Y$ must be in $F.$}.

		Next we formalize  what we mean by  ``performing the linear change of variables $x \leftarrow Mx$ in \eqref{problem-sdp}" for some invertible matrix $M.$ 
		\begin{Definition} \label{definition-reform}  We say that we  {\em reformulate } \eqref{problem-sdp} 
			if we apply to it some of the following operations (in any order): 
			\benum
			\item \label{exch} Exchange $A_i$ and $A_j, \,$ where $i$ and $j$ are distinct indices in $\{1, \dots, m \}.$ 
			\item \label{trans} Replace $A_i$ by $\lambda A_i  + \mu A_j, \,$ where $i$ and $j$ are distinct indices in $\{1, \dots, m \}, \,$
			$\lambda$ and $\mu$ are reals, and $\lambda \neq 0.$ 
			\item \label{rotate} Replace all $A_i$ by 
			$T^{\top} A_iT  \, $ and $B$ by $T^{\top} B T, \, $ where $T$ is a suitably chosen invertible matrix. 
			\eenum
			We also say that by reformulating  \eqref{problem-sdp}  we obtain a {\em reformulation}.
			\end{Definition}
Reformulations were  originally introduced to study various pathologies in SDPs, for example, 
unattained optimal values and duality gaps \cite{pataki2019characterizing}; and infeasibility \cite{LiuPataki:15}. 
In this work we show that they help understand another classical pathology, 
exponential size solutions.

We next clarify some technicalities about reformulations. First, 
		the above definition of a reformulation slightly differs from the one 
		in \cite{pataki2019characterizing},  where we also permit replacing $B$ by $B + \lambda A_i$ for some $i$ index and $\lambda$ real number.
			Second, operations \eqref{exch} and \eqref{trans} can be viewed  as elementary row operations on a dual type system, say, on  
			$$
			A_i \bullet Y = 0 \,\, \text{for} \,\, i=1, \dots, m.
			$$
			Third, operation \eqref{rotate} does not influence the values of the $x_i,$ since $x \in \rad{m}$ is feasible in \eqref{problem-sdp} before 
		we apply operation 
		\eqref{rotate},  if and only if it is feasible in it afterwards.  Thus, we only use
		operation \eqref{rotate}  to put \eqref{problem-sdp} into a more convenient looking form.

		We will rely on the following statement about the connection of $\psdn$ and a linear subspace.  
		 It is a special case of a classical, more general statement about the intersection of a linear subspace and a  convex cone: see e.g., 
	  \cite[Theorem 2]{LuoSturmZhang:97}.  
		\begin{Proposition} \label{prop-gordan-stiemke} 
			Suppose $L$ is a linear subspace of $\symn.$ Then exactly one of the following two alternatives is true:
			\begin{enumerate}
				\item \label{alt-psd-nonzero} There is a nonzero positive semidefinite matrix in $L.$ 
					\item \label{alt-pd} There is a positive definite matrix in $L^\perp.$ 
			\end{enumerate}
		\qed
		\end{Proposition}

\section{Main results and proofs} \label{section-main}

\subsection{Reformulating \eqref{problem-sdp} and statement of Theorem \ref{thm-main}} 
\label{subsection-state-thm1}

In our first lemma we present an algorithm to reformulate \eqref{problem-sdp} into a more convenient looking form. 
The algorithm is a simplified version of the algorithm in \cite{LiuPataki:15}  \footnote{
	The algorithm of Lemma \ref{lemma-reformulation} uses only Proposition \ref{prop-gordan-stiemke}, whereas the algorithm of
\cite{LiuPataki:15} relies on a more involved  theorem of the alternative.}. 
Both algorithms are specialized facial reduction algorithms
 applied to the dual semidefinite system defined in \eqref{eqn-A_iY=0}.

\begin{Lemma} \label{lemma-reformulation} 
	The problem 
	\eqref{problem-sdp} has a reformulation 
\begin{equation} \label{problem-sdpprime} \tag{\mbox{$\mathit{P'}$}} 
x_1 A_1' + \dots + x_k A_k' + x_{k+1} A_{k+1}' + \dots + x_m A_m' + B' \succeq 0 
\end{equation}
with the following properties:
\begin{itemize}
	\item $k$ is a nonnegative integer, and $(A_1', \dots, A_k')$ is a regular facial reduction sequence.
	\item If $r_1, \dots, r_k$ is the size of the 
	identity block in $A_1', \dots, A_k', $ respectively, then 
	$n - r_1 - \dots - r_k$ is the maximum rank of a  matrix in 
	\begin{equation} \label{eqn-A_iY=0} 
		\{  \, Y \succeq 0 \, | \, A_i \bullet Y = 0 \; \text{for} \; i=1, \dots, m \,   \}.
	\end{equation}
\end{itemize}
\end{Lemma}

\pf{} 
We will reformulate \eqref{problem-sdp} in several steps. 
To start, we let $L$ be the linear span of $A_1, \dots, A_m$ and apply Proposition \ref{prop-gordan-stiemke}.
If  item \eqref{alt-pd} holds,  we  let $k=0, \, A_i' = A_i$ for all $i, \, $ $B' = B, \, $  and stop.

If item \eqref{alt-psd-nonzero} holds, then we choose a nonzero  psd matrix $V = \sum_{i=1}^m \lambda_i A_i$  in 
$L$ and assume $\lambda_1 \neq 0$ without loss of generality.  
 We then choose a $T$ invertible matrix so that 
 \begin{equation} \label{eqn-TVT} 
T^\top V T = \bpx I_{r_1} & 0 \\
                                    0     & 0 \epx,
\end{equation}
where $r_1$  is the rank of $V.$ 
We let $A_1' := T^\top V T, A_i' := T^\top A_i T$ for $i \geq 2,$ and $B' = T^\top B T.$

Let $r$ be the maximum rank of a psd matrix in $L^\perp$ (i.e., in the set defined in \eqref{eqn-A_iY=0}). Also, let $L_{\mathrm new}$ be the linear span of $A_1', \dots, A_m'.$ We claim that 
$r$ is also the maximum rank of a psd matrix in  $L_{\mathit new}^\perp.$ 
For that,  let us choose  a rank $r$ matrix, say $Y, $ in $L^\perp \cap \psd{n}.$ 
Then for all $i$ we have 
$$
0 \, = \, A_i \bullet Y \, = \, T^\top A_i T \bullet T^{-1} Y T^{- \top}, 
$$
where the last equality is from the definition of the $\bullet$ product and the properties of the trace.
Thus $T^{-1} Y T^{- \top}$ is in  $L_{\mathit new}^\perp$ and has rank $r.$ 
Similarly, from any psd matrix in $L_{\mathrm new}^\perp$ we can construct a psd matrix in $L^\perp$ with the same rank. This proves our claim.

Suppose that $Y \in L_{\mathrm new}^\perp \cap \psdn; \, $ then $A_1' \bullet Y = 0.$ Since $A_1'$ is now the $T^\top V T$ matrix given in 
\eqref{eqn-TVT}, the sum of the first $r_1$ diagonal elements of $Y$ is zero. Since $Y$ is psd, the first $r_1$ rows and columns of $Y$ are zero.

We next construct an SDP 
$$
\sum_{i=2}^m x_i F_i + G \succeq 0,
$$
where $F_i$ is obtained from $A_i'$ by deleting the first $r_1$ rows and columns for $i=2, \dots, m, \,$ 
and $G$ is obtained from $B'$ in the same manner. 
By the above argument the maximum rank of a matrix in $\{Z \succeq 0 : F_i \bullet Z = 0 \, (i=2, \dots, m) \}$
is also $r,$ so we can  proceed in a similar manner with this smaller SDP.  When our process stops, 
we  have the required reformulation.
\qed

The reader may wonder why we require $(A_1', \dots, A_k')$ to be a regular facial reduction sequence in 
	\eqref{problem-sdpprime}. Will we use them to verify that any $Y$ in the set 
	\begin{equation} \label{eqn-A_iprimeY=0} 
	\{  \, Y \succeq 0 \, | \, A_i '\bullet Y = 0 \; \text{for} \; i=1, \dots, m \,   \}
\end{equation}
has its first $n- r_1 - \dots - r_k$ rows and columns equal to zero? We could  indeed use them for this purpose, 
by an argument similar to the one  after Definition \ref{definition-regfr} 
(the $A_i'$  would play the role of the $C_i$). 
However, interestingly, such a $Y$ will never appear in our arguments in Lemmas \ref{lemma-xk-arb-large},  \ref{lemma-tailindex}, and later. 
Instead, we will use the staircase structure of the $A_i'$ to prove results about large size solutions in \eqref{problem-sdpprime}.

From now on we assume that 
\begin{center} 
	\framebox[4.5in]{$k$  is the smallest integer that satisfies the requirements of Lemma 
		\ref{lemma-reformulation}.} 
\end{center}

 Using the terminology of facial reduction, $k$ is the {\em singularity degree of} the dual system 
 \eqref{eqn-A_iY=0} \footnote{We can also define the singularity degree of \eqref{problem-sdp}. Since this problem is strictly feasible, its singularity degree is  just zero.}. This concept was originally introduced by Sturm in \cite{Sturm:00} 
 and used to derive error bounds, namely, bounds on the distance of a point from the feasible set of an SDP.  
  For a broad 
 generalization of Sturm's result to conic systems over 
 so-called amenable cones, see a recent result by  Louren{\c{c}}o \cite{lourencco2017amenable}.

Note  that since \eqref{problem-sdp} is strictly feasible, so is \eqref{problem-sdpprime}.
Since we will focus on the leading $k$ variables in \eqref{problem-sdpprime}, 
for the rest of the paper we fix $(\bar{x}_{k+1}, \dots, \bar{x}_m)$  such that 
\begin{center} 
	\framebox[4.5in]{$(x_1, \dots, x_k, \bar{x}_{k+1}, \dots, \bar{x}_m)$ is strictly feasible in \eqref{problem-sdpprime} for some $x_1, \dots, x_k.$}
\end{center}   

 From now on we will  say that a number is a {\em constant}, if it depends only on 
  the $\bar{x}_i, \, $ the $A_i, \,$ and $B.$ Theorem \ref{thm-main} will rely on such constants. 
 	
We now formally state our main result.

\begin{Theorem}   \label{thm-main}     
	Let \eqref{problem-sdpprime} be the reformulation of \eqref{problem-sdp} obtained in Lemma \ref{lemma-reformulation}, $k$ the singularity 
	degree of the dual system \eqref{eqn-A_iY=0}, and assume $k \geq 2.$  Then 
	
	\begin{enumerate}
		\item 	 \label{thm-main-1} There is $(x_1, \dots, x_k)$ such that $(x_1, \dots, x_k, \bar{x}_{k+1}, \dots, \bar{x}_m)$ is strictly feasible in 
		\eqref{problem-sdpprime} and  $x_k$ is arbitrarily large.
		\item  \label{thm-main-2}  If  $(x_1, \dots, x_k, \bar{x}_{k+1}, \dots, \bar{x}_m)$ is strictly feasible in \eqref{problem-sdpprime} and 
			$x_k$ is sufficiently large, then 
		\begin{equation} \label{eqn-xj-xj+1}
		\ba{rcl} 
		x_j & \geq & d_{j+1} x_{j+1}^{\alpha_{j+1}} \, \text{for \,} j=1, \dots, k-1,
		\ena 
		\end{equation}
		where 
		\begin{equation} \label{eqn-alphaj-bounds} 
		2  \geq \alpha_{j+1}  \geq  1 + \dfrac{1}{k-j} \,\,\,  \text{for} \, j=1, \dots, k-1.  
		\end{equation}
		Here the $d_j$ and $\alpha_j$ are positive constants. 
	
\end{enumerate} 
\qed
\end{Theorem}

Note that even $k \geq 1$ easily  implies that the feasible set of \eqref{problem-sdpprime} (or equivalently, that of \eqref{problem-sdp}) 
	is unbounded. Indeed, since $A_1' \succeq 0, \, $ we can add an arbitrarily large multiple of the first unit vector to any feasible solution of \eqref{problem-sdpprime}, and stay feasible.
 Of course, Theorem \ref{thm-main} proves much more than the feasible set of 
		\eqref{problem-sdpprime} is unbounded: it proves a hierarchy among the variables in \eqref{problem-sdpprime}.

The proof  of  Theorem \ref{thm-main} has  three main  parts. First, in Lemma 
\ref{lemma-xk-arb-large} we prove item \eqref{thm-main-1}, that in strictly feasible solutions of  \eqref{problem-sdpprime} 
we can have arbitrarily large $x_k.$ 
 Lemma \ref{lemma-tailindex} is a technical statement about  a certain parameter, called the 
{\em tail-index} of the $A_i';$ this parameter depends on where the  nonzero blocks of the $A_i'$ are. 

\co{Lemma \ref{lemma-tailindex} is a technical statement about  nonzero blocks in the $A_i', \,$ which yields a lower bound on a certain parameter, called their 
{\em tail-index}.}  

 In the second part, Lemma \ref{lemma-derive-poly}  deduces from \eqref{problem-sdpprime} a 
set of 
quadratic inequalities. These are typically ``messy", namely they look like  
$$
(x_1 + x_2 + x_3)(x_4 + 10 x_5) >  (x_2 - 3 x_4)^2.
$$
Third, in Lemma \ref{lemma-xj>=dj+1} from these messy inequalities we first derive ``cleaned up" versions, such as 
$$
x_1 x_4  >  \const x_2^2.
$$
Then from these cleaned up inequalities 
we deduce the inequalities \eqref{eqn-xj-xj+1} and a recursive formula to compute  the $\alpha_j.$ 
Next,  Lemma \ref{lemma-alphaj-monotone} proves that the $\alpha_j$ exponents are a monotone function 
 of the tail-indices of the $A_j'.$ Finally, Lemma \ref{lemma-alphaj-smallest} 
 shows that minimal tail-indices of the $A_j'$ give the smallest possible $\alpha_j$  exponents.  
 We then combine all lemmas and prove  Theorem \ref{thm-main}.

Before we  get to the proof, we illustrate Theorem \ref{thm-main} 
via  two extreme examples.  Recall that  Theorem \ref{thm-main} is about strictly feasible solutions of \eqref{problem-sdpprime},  
in which $x_k$ must be large enough.
However, the following SDP examples are fairly simple, and in all of them we can derive interesting quadratic inequalities that hold for all feasible solutions.  
  
\begin{Example} (Khachiyan SDP) \label{example-khanch-sdp} 
	Consider the SDP 
	\begin{equation} \label{problem-khach-sdp} 
	\tag{\mbox{${\mathit  Kh \mhyphen SDP}$}}
	\begin{pmatrix}
	\mrbmath{x_1}  	 &           &            &     & \mrbmath{x_2}  \\
	& \mbbmath{x_2}   &           &     & \mbbmath{x_3}  \\
	&           & \mgbmath{x_3}    &     & \mgbmath{x_4}  \\  
	&           &            &  {x_4}   &         \\
	\mrbmath{x_2}   & \mbbmath{x_3}  & \mgbmath{x_4}         &   	& {1}  
	\end{pmatrix} \succeq 0, 
	\end{equation}
which can be written  in the form of \eqref{problem-sdpprime} with the $A_i'$ matrices 
given below and $B'$ the matrix whose lower right corner is $1$ and the remaining elements are zero:
\begin{equation} \nonumber 
\underbrace{\bpx 1 & & & & \\
	& 0 & & & \\
	& & 0 & & \\
	& & & 0 & \\
	& & & & 0  
	\epx}_{A_1'},  \underbrace{\bpx  0 & & & & 1 \\
	& 1 & & & \\
	& &0 & & \\
	& & & 0 & \\
	1 & & & & 0
	\epx}_{A_2'},   \underbrace{\bpx  0 & & & &  \\
	& 0 & & & 1 \\
	& & 1 & & \\
	& & &0 & \\
	& 1 & & & 0
	\epx}_{A_3'},  \underbrace{\bpx 0 & & & & \\
	& 0 & & & \\
	& & 0 & & 1\\
	& & & 1 & \\
	& & 1 & & 0  
	\epx}_{A_4'}. \co{,   \underbrace{\bpx 0 & & & & \\
		& 0 & & & \\
		& & 0 & & \\
		& & & 0 & \\
		& & & & 1  
		\epx}_{B'},}
\end{equation} 
The subdeterminants in \eqref{problem-khach-sdp} with three red, three blue, and three green corners, respectively,
	give the inequalities   
	\begin{equation} \label{eqn-khach-quadratic} 
	x_1 \geq x_2^2, \, x_2 \geq x_3^2, \, x_3 \geq x_4^2
	\end{equation}
	that appear in \eqref{problem-khachiyan}. So the exponents in the inequalities \eqref{eqn-khach-quadratic} 
		are the largest permitted by our bounds \eqref{eqn-alphaj-bounds}. 
		
	(For simplicity we constructed this SDP, so its feasible set does not imply the inequality $x_4 \geq 2,$ which does appear in \eqref{problem-khachiyan}.)

\end{Example}
 What is $k$ in Example \ref{example-khanch-sdp}? By definition, it is the singularity degree of 
	\begin{equation} \label{elaine} 
		\{ Y \succeq 0: A_i' \bullet Y = 0 \; {\rm for} \; i=1,\dots, 4 \, \}
	\end{equation} 
and we will next show that $k=4.$ For that, we first observe that the  maximum rank of a matrix in this system is $1.$ 
Then we consider a reformulation of \eqref{problem-khach-sdp} 
\begin{equation} \label{eqn-sum-barAi} 
\sum_{i=1}^4 x_i \bar{A}_i + \bar{B} \succeq 0,
\end{equation}
with two properties. First, for some $k \leq 4$ the sequence 
$(\bar{A}_1, \dots, \bar{A}_k)$ is a regular facial reduction sequence. Second, the sizes 
of the identity blocks in the $\bar{A}_i$ sum to $5-1 = 4$ (note that now $n=5$). 
Since $k$ is the singularity degree of \eqref{eqn-sum-barAi}, it is minimal, so the identity blocks in   $(\bar{A}_1, \dots, \bar{A}_k)$ 
are nonempty. 
Thus $\bar{A}_1$ is a positive multiple of $A_1',$ since $A_1'$ is the only nonzero psd matrix in the linear span of the $A_i'. \, $ 
Similarly, it follows by induction for $i=1, \dots, k$ that each 
 $\bar{A}_i$ is a positive multiple of $A_i'$ plus a linear combination of $A_{i-1}', \dots, A_1'.$ Thus $k=4$ follows.

We note in passing that the feasible sets  of \eqref{problem-khach-sdp} and of 
the derived quadratic inequalities \eqref{eqn-khach-quadratic} are not equal. For example  $x=(256,16,4,2)$ is not feasible in 
\eqref{problem-khach-sdp},  but  is feasible in  \eqref{eqn-khach-quadratic}. 
However, we can  construct an SDP that exactly represents \eqref{problem-khachiyan}, as we described in Section \ref{sect-intro}. In that SDP a  straightforward argument like the one we gave above proves
that $k=m.$ 

We next discuss whether we need to assume that a strictly feasible solution exists in \eqref{problem-sdp}, in order to derive Theorem 
\ref{thm-main}. On one hand, there are semidefinite programs 
which have  no strictly feasible solutions, nor do they exhibit the hierarchy among the leading variables 
seen in  Theorem \ref{thm-main}.  
Indeed, we obtain such an SDP if in  \eqref{problem-khach-sdp} we change $x_1$ to $x_1 + 1$ and the $1$ entry in the bottom right corner to $0$.  This new SDP is represented by the same $A_i'$ matrices, so the associated singularity degree of  \eqref{elaine} is still four. 
Further, this new SDP is not strictly feasible, and $x_2 = x_3 = x_4 = 0$ holds in any feasible solution, but $x_1$ can be $-1.$ 
{We can similarly create such an SDP with an arbitrary number of variables.}

On the other hand, there are SDPs with no strictly feasible solution, which, however, do 
have a Khachiyan type hierarchy among the variables (and hence 
large size solutions). For that, we only need to take  an SDP 
with a  Khachiyan type hierarchy, 
and simply 
add 
all-zero rows and columns.  

	\begin{Example}          (Mild SDP)  \label{example-mild} 
		As a counterpoint to \eqref{problem-khach-sdp} we next consider a mild SDP (we will see soon why we call it ``mild")
			\begin{equation} \label{problem-mild-sdp} \tag{\mbox{${\mathit  Mild \mhyphen SDP}$}}
			\begin{pmatrix}
			\mrbmath{x_1} 	 &           & \mrbmath{x_2}           &     &   \\
			& \mbbmath{x_2}   &           &  \mbbmath{x_3}   &   \\
			\mrbmath{x_2}     &           & \mgbmath{x_3}     &     & \mgbmath{x_4}  \\
			& \mbbmath{x_3}         &            & \mblackbmath{x_4}   &         \\
			&   & \mgbmath{x_4}          &  & 	{\bf 1}       
			\end{pmatrix} \succeq 0. 
			\end{equation}
				We naturally write \eqref{problem-mild-sdp} in the form of \eqref{problem-sdpprime} with the 
				$A_i'$ matrices shown below and 
				$B'$ the matrix whose lower right corner is $1$ and the remaining elements are zero:
				\begin{equation} \nonumber 
				\underbrace{\bpx 1 & & & & \\
					& 0 & & & \\
					& & 0 & & \\
					& & & 0 & \\
					& & & & 0  
					\epx}_{A_1'},  
				\underbrace{\bpx  0 &  & 1 & &   \\
					& 1 & & & \\
					1 & &0 & & \\
					& & & 0 & \\
					& & & & 0
					\epx}_{A_2'},   \underbrace{\bpx  0 & & & &  \\
					& 0 &  & 1 &  \\
					& & 1 & & \\
					& 1 & &0 & \\
					&  & & & 0
					\epx}_{A_3'},  \underbrace{\bpx 0 & & & & \\
					& 0 & & & \\
					& & 0 & & 1\\
					& & & 1 & \\
					& & 1 & & 0  
					\epx}_{A_4'}. 
				\end{equation}
				In \eqref{problem-mild-sdp} the subdeterminants with three red, three blue, and three green corners, respectively, yield 
			the inequalities 
			\begin{equation} \label{eqn-mild-quadratic} 
			x_1 x_3  \geq x_2^2, \, x_2 x_4 \geq  x_3^2, \, x_3  \geq  x_4^2.  
			\end{equation}
					Next from the inequalities in \eqref{eqn-mild-quadratic} we derive
				\begin{equation} \label{eqn-mild-quadrati-corollary} 
				x_1  \geq x_2^{4/3}, \, x_2 \geq x_3^{3/2}, \, x_3  \geq  x_4^2
				\end{equation}
				as follows. We first copy the last  inequality $x_3 \geq x_4^2$  from \eqref{eqn-mild-quadratic} to \eqref{eqn-mild-quadrati-corollary}.
				Next we plug $x_3^{1/2} \geq x_4$ into the middle inequality in \eqref{eqn-mild-quadratic} to get  $x_2 \geq x_3^{3/2}.$ 
				We finally  raise both sides of this last  inequality to the power of $2/3$ and plug it into the first inequality in 
				\eqref{eqn-mild-quadratic} to  deduce $x_1 \geq x_2^{4/3}.$ 
			
			To summarize,  the exponents in the derived inequalities \eqref{eqn-mild-quadrati-corollary} 
			are the smallest permitted by our bounds \eqref{eqn-alphaj-bounds}. 
			
	\end{Example} 
We invite the reader to verify that $k=4$ holds in Example \ref{example-mild}; this can be done just like in Example \ref{example-khanch-sdp}.

	To illustrate the difference between \eqref{problem-khachiyan} and the inequalities derived from 
	\eqref{problem-mild-sdp}, we show the set defined by the inequalities 
	\begin{equation} \label{eqn-mild-quadratic-for-pic} 
	x_1 x_3  \geq x_2^2, \, x_2  \geq  x_3^2, \, 2 \geq x_3  \geq  0
	\end{equation}
	  on the right in Figure \ref{figure-together}.  Note that the set defined by 
	  \eqref{eqn-mild-quadratic-for-pic} is a three dimensional version of the set given 
	  in \eqref{eqn-mild-quadratic}, which we normalized by adding upper and lower bounds on 
	  $x_3.$

\subsection{Proof of Theorem \ref{thm-main}}
\label{subsection-proofs}

In Lemmas \ref{lemma-xk-arb-large}--\ref{lemma-derive-poly} we will use the following notation:
\begin{equation} \label{eqn-def-r_i-I_i}
\ba{rcl} 
r_j & = &  \text{size of the identity block in } \, A_j' \, \text{ for} \,   j=1, \dots, k, \\
\I_1 & := & \{ 1, \dots, r_1 \},  \\
\I_2 & := &  \{ r_1 + 1, \dots, r_1 + r_2 \},  \\
& \vdots & \\
\I_k & := & \{ r_1 + \dots + r_{k-1} +1, \dots, r_1 + \dots + r_k    \}, \\ 
\I_{k+1} & := &  \{ r_1 + \dots + r_k +1, \dots, n \}.
\ena 
\end{equation}

 \begin{Lemma} \label{lemma-xk-arb-large} 
 There is $(x_1, \dots, x_k)$ such that $(x_1, \dots, x_k, \bar{x}_{k+1}, \dots, \bar{x}_m)$ is strictly feasible in 
 	\eqref{problem-sdpprime} and $x_k$ is arbitrarily large.  
 \end{Lemma}
 \pf{} 
 Let 
 $$
 Z := \sum_{i=k+1}^m \bar{x}_i A_i' + B'.
 $$
 Since 
 there is $x_1, \dots, x_k$ such that $\sum_{i=1}^k x_i A_i' + Z \succ 0, \,$ and 
 $A_i'(\I_{k+1}) =0$ for $i=1, \dots, k$ 
 we see that 
 $$Z(\I_{k+1}) \succ 0.$$ 
Recall that $A_k' (\I_k) = I$ and the other elements of $A_k'(\I_k \cup \I_{k+1})$ are zero. Hence 
 by the definition of positive definiteness (a symmetric matrix $G$ is positive definite 
 if $x^\top G x > 0$ for all nonzero $x$) we see that 
 the $\I_k \cup \I_{k+1}$ diagonal block of $x_k A_k' + Z$ is positive definite when $x_k$ is large enough. 
 Similarly, for any such $x_k$ there is  $x_{k-1}$ so the $\I_{k-1} \cup \I_k \cup \I_{k+1}$ diagonal block of 
 $x_{k-1} A_{k-1}' + x_k A_k' + Z$ is positive definite. We construct $x_{k-2}, \dots, x_1$ in a similar manner. 
 \qed
 
	The proof of Lemma \ref{lemma-xk-arb-large} is partly inspired by the paper of Louren{\c{c}}o et al \cite{lourencco2016structural}, which used a similar process to construct a nearly feasible solution to a weakly infeasible semidefinite program. 

The  proof of Lemma \ref{lemma-xk-arb-large} also partially answers the representation 
question \eqref{question-exp-represent}.  To explain how,  for the moment let us ignore the requirement that we need to choose $x_k$ to be large and just focus
on completing $(\bar{x}_{k+1}, \dots, \bar{x}_m)$ to a strictly feasible solution. 
The proof that the required $(x_1, \dots, x_k)$ {\em could be} computed is fairly simple, and it is 
illustrated on Figure \ref{figure-add-xk}, where the red blocks stand for the larger and larger blocks that we make positive definite.  So we can convince ourselves that $(x_1, \dots, x_k)$ 
{\em exist}, even without computing their   actual values. 
\begin{figure}[H] 
	\begin{center} 
		\begin{tikzpicture}
\tikzstyle{every left delimiter}=[xshift=0.9ex]
\tikzstyle{every right delimiter}=[xshift=-0.9ex]

\newcommand{\ph}[1]{\phantom{#1}}
\usetikzlibrary{shapes.misc,shadows} 
\usetikzlibrary{decorations.pathreplacing,angles,quotes}
\usetikzlibrary{matrix}
\usetikzlibrary{backgrounds, fit}

\matrix[
    matrix of math nodes,
    row sep=0.5ex,
    column sep=0.5ex,
    left delimiter=(,right delimiter=),
    nodes={text width=0.8em, text height=0.8em, text depth=0.5ex, align=center},
    every node/.style={scale=0.9}
    ] (m) at (0,0)
    {
    $\times$ & $\times$ & $\times$ & $\times$ \\
    $\times$ & $\times$ & $\times$ & $\times$ \\
    $\times$ & $\times$ & $\times$ & $\times$ \\
    $\times$ & $\times$ & $\times$ & $+$\\
    };
    \begin{scope}[on background layer]
        \filldraw[red!70] (m-4-4.north west) -- 
        (m-4-4.north east) -- (m-4-4.south east)-- (m-4-4.south west) --
        cycle;
    \end{scope} 
    
\matrix[
    matrix of math nodes,
    row sep=0.5ex,
    column sep=0.5ex,
    left delimiter=(,right delimiter=),
    nodes={text width=0.8em, text height=0.8em, text depth=0.5ex, align=center},
    every node/.style={scale=0.9}
    ] (m) at (4.4,0)
    {
    $\times$ & $\times$ & $\times$ & $\times$ \\
    $\times$ & $\times$ & $\times$ & $\times$ \\
    $\times$ & $\times$ & $+$ & $\times$ \\
    $\times$ & $\times$ & $\times$ & $+$\\
    };
    \begin{scope}[on background layer]
        \filldraw[red!70] (m-3-3.north west) -- 
        (m-3-4.north east) -- (m-4-4.south east)-- (m-4-3.south west) --
        cycle;
    \end{scope}
    
\matrix[
    matrix of math nodes,
    row sep=0.5ex,
    column sep=0.5ex,
    left delimiter=(,right delimiter=),
    nodes={text width=0.8em, text height=0.8em, text depth=0.5ex, align=center},
    every node/.style={scale=0.9}
    ] (m) at (4.4*2,0)
    {
    $\times$ & $\times$ & $\times$ & $\times$ \\
    $\times$ & $+$ & $\times$ & $\times$ \\
    $\times$ & $\times$ & $+$ & $\times$ \\
    $\times$ & $\times$ & $\times$ & $+$\\
    };
    \begin{scope}[on background layer]
        \filldraw[red!70] (m-2-2.north west) -- 
        (m-2-4.north east) -- (m-4-4.south east)-- (m-4-2.south west) --
        cycle;
    \end{scope}

    \draw[decoration={brace,raise=2.5pt},decorate] (0.5,1.05) -- node[above=5pt] {\footnotesize$\mathcal{I}_{k+1}$} (1.1,1.05);
    \draw[decoration={brace,raise=2.5pt,mirror},decorate] (0.5,-1.05) -- node[below=5pt] {\footnotesize$\succ 0$} (1.1,-1.05);
    \draw[decoration={brace,raise=2.5pt,mirror},decorate] (-1.1,-1.8) -- node[below=5pt] {$Z$} (1.1,-1.8);
    
    \node[] at (2.2,0) {$\longmapsto$};
    \node[] at (2.2,0.35) {\scriptsize$+x_kA_k'$};
    \node[] at (2.2,-0.35) {\scriptsize$x_k\gg 0$};
    
    \draw[decoration={brace,raise=2.5pt},decorate] (4.4+0.5-0.525,1.05) -- node[above=5pt] {\footnotesize$\mathcal{I}_{k}$} (4.4+1.1-0.525,1.05);
    \draw[decoration={brace,raise=2.5pt,mirror},decorate] (4.4+0.5-0.525,-1.05) -- node[below=5pt] {\footnotesize$\succ 0$} (4.4+1.1,-1.05);
    \draw[decoration={brace,raise=2.5pt,mirror},decorate] (4.4+-1.1,-1.8) -- node[below=5pt] {$x_kA_k'+Z$} (4.4+1.1,-1.8);
    
    \node[] at (4.4+2.2,0) {$\longmapsto$};
    \node[] at (4.4+2.2,0.35) {\scriptsize$+x_{k-1}A_{k-1}'$};
    \node[] at (4.4+2.2,-0.35) {\scriptsize$x_{k-1}\gg 0$};
    
    \draw[decoration={brace,raise=2.5pt},decorate] (4.4*2+0.5-0.525*2,1.05) -- node[above=5pt] {\footnotesize$\mathcal{I}_{k-1}$} (4.4*2+1.1-0.525*2,1.05);
    \draw[decoration={brace,raise=2.5pt,mirror},decorate] (4.4*2+0.5-0.525*2,-1.05) -- node[below=5pt] {\footnotesize$\succ 0$} (4.4*2+1.1,-1.05);
    \draw[decoration={brace,raise=2.5pt,mirror},decorate] (4.4*2+-1.1,-1.8) -- node[below=5pt] {$x_{k-1}A_{k-1}'+x_kA_k'+Z$} (4.4*2+1.1,-1.8);
    
    \node[] at (4.4*2+2.2,0) {$\longmapsto$};
    \node[] at (4.4*2+2.2,0.35) {\scriptsize$+x_{k-2}A_{k-2}'$};
    \node[] at (4.4*2+2.2,-0.35) {\scriptsize$x_{k-2}\gg 0$};
    
    \node[] at (4.4*2+2.2+1.2,0) {$\dots$};
    
\end{tikzpicture}
	\end{center} 
		\caption{Verifying that $x_1, \dots, x_k$ exist, without computing them}
		\label{figure-add-xk}
\end{figure}
From now on we will assume 
 \begin{equation} \label{eqn-sum-rj} 
 r_1 + \dots + r_k < n,
 \end{equation}
 and we claim that we can do so  without loss of generality. 
 Indeed, suppose that the sum of all the $r_j$ is $n.$ Then an argument like in the proof of Lemma \ref{lemma-xk-arb-large} proves that 
 $A_1', \dots, A_k'$ have a positive definite linear combination. Hence the singularity degree of 
 \eqref{eqn-A_iY=0} is actually just $1; \, $ but we assumed $k \geq 2.$ 

By \eqref{eqn-sum-rj} we see that  $\I_{k+1} \neq \emptyset.$ 

To motivate our next definition we compare our two extreme examples from two viewpoints. 
From the first viewpoint we see that in   \eqref{problem-khach-sdp} the $x_j$ variables in the upper offdiagonal positions are more to the right than in \eqref{problem-mild-sdp}. 
From the second viewpoint, in the  inequalities \eqref{eqn-khach-quadratic} 
derived from   \eqref{problem-khach-sdp}  the exponents are larger than in the inequalities \eqref{eqn-mild-quadrati-corollary} derived from \eqref{problem-mild-sdp}. 
We will see that these two facts are closely connected, so in the next definition we capture ``how far to the right the $x_j$ are in upper 
offdiagonal positions." 
  
\begin{Definition} \label{definition-tailindex} 
We define the {\em tail-index} $t_{j+1}$ of $A_{j+1}'$ for $j=1, \dots, k-1$ as 
	\begin{equation} \label{eqn-tail-index} 
	t_{j+1} \, :=  \, \max \, \{ \, t  \, : \, A_{j+1}'(\I_j, \I_t)    \neq 0  \}.  \; 
	\end{equation}
\end{Definition} 
In words, $t_{j+1}$ is the index of the rightmost nonzero block of columns ``directly above" the identity block in $A_{j+1}'.$ We illustrate the tail-index on Figure \ref{figure-tailindex}. Here and in later figures the $\bullet$ blocks are nonzero, and we separate the columns indexed by $\I_{k+1}$ from the other columns by double vertical lines.

We further illustrate Definition \ref{definition-tailindex} using Examples \ref{example-khanch-sdp} and \ref{example-mild}. In both of these, we have $ \I_j = \{ j \}$ for $j=1,\dots, 5.$ 
Thus, the tail-indices are $t_2=t_3=t_4=5 \, $ in \eqref{problem-khach-sdp},  whereas they are 
	$t_2 = 3, \, t_3 = 4, \, t_4 = 5$ in \eqref{problem-mild-sdp}.

	\begin{figure}
		\begin{small}   
		$$
		A_{j+1}' \, = \, \begin{pmatrix}[c|c|c|c|c|c||c]
			{\mbox{$\,\,\,\,\ti \,\,\,\,$}} 	& \bovermat{$\I_{j}$}{\mbox{$\,\,\, \ti \,\,\,\,$}}	& \bovermat{$\I_{j+1}$}{\mbox{$\,\,\, \ti \,\,\,\,$}} & {\mbox{$\,\,\,\, \ti \,\,\,$}}	& \bovermat{$  \I_{t_{j+1}}$}{\mbox{$\,\,\,\,\ti \,\,\,$}}	 & {\mbox{$\,\,\,\,\ti\,\,\,\,$}} & \bovermat{$\I_{k+1}$}{\mbox{$\,\,\,\,\ti \,\,\,$}}        \\ \hline 
			\,\,\,\ti \,\, & \ti & \ti & \ti & \bullet  &   &  \\ \hline 
			\,\,\,\ti \,\, & \ti & I &   &    &  &   \\ \hline 
			\,\,\,\ti \,\, & \ti &  &  &    &  &  \\ \hline 
			\ti & \bullet &  &  & & &  \\ \hline  
			\ti &   &  &  & & &  \\ \hline \hline 
			\ti &  &  & & & & 
			\end{pmatrix}. 
		$$
		\end{small} 
		\caption{The tail-index of $A_{j+1}'$} 
		\label{figure-tailindex} 
	\end{figure}

	\begin{Lemma}  \label{lemma-tailindex}
		$$
		t_{j+1} > j+1 \,\, \text{for} \,\, j=1, \dots, k-1.
		$$
	\end{Lemma}
\pf{} 
	We will use the following notation: 
	for $r, s  \in \{1, \dots, k+1 \}$ such that $r \leq s$ we let 
		\begin{equation} \label{eqn-define-I_i_j} 
			\I_{r:s}  :=  \I_r \cup \dots \cup \I_s. 
		\end{equation}
		Let $j \in \{1, \dots, k-1\}$ be arbitrary. 
		To help with  the proof, we picture $A_{j}'$ and $A_{j+1}'$ in equation \eqref{bla}.
	As always, 
	the empty blocks are zero, and the $\ti$ blocks are arbitrary. The blocks marked by $\otimes$ are 
	$A_{j+1}'( \I_{j}, \I_{(j+2):(k+1)})$ and its symmetric counterpart.   We will prove that these blocks are nonzero and this will prove  our lemma. 
	
		\vspace{.1cm}      
			\beq \label{bla}    
		A_{j}' \, = \, \begin{pmatrix}[c|c|c|c]
			\bovermat{$\I_{1:(j-1)}$}{\mbox{$\,\,\,\, \times \,\,\,\,$}} 	& \bovermat{$\I_{j}$}{\mbox{$\,\,\,\times\,\,\,$}}	& \bovermat{$ \I_{j+1}$}{\mbox{$\,\,\,\times\,\,\,$}}	& \bovermat{$  \I_{(j+2):(k+1)}$}{\mbox{$\,\,\,\,\quad\times\quad\,\,\,\,\,$}}	\\ \hline 
			\ti & I &  &    \\ \hline
			\ti &  & &    \\ \hline
			\ti & \pha{0} &  &  
		\end{pmatrix}, \,  
		A_{j+1}'  \, = \, \begin{pmatrix}[c|c|c|c]
			\bovermat{$\I_{1:(j-1)}$}{\mbox{$\,\,\,\, \times \,\,\,\,$}} 	& \bovermat{$\I_{j}$}{\mbox{$\,\,\, \times\,\,\,$}}	& \bovermat{$ \I_{j+1}$}{\mbox{$\,\,\,\times\,\,\,$}}	& \bovermat{$\I_{(j+2):(k+1)}$ }{\mbox{$\,\,\,\,\quad\times\quad\,\,\,\,\,$}}	\\ \hline 
			\ti & \ti & \ti & \otimes  \\ \hline
			\ti & \ti & I &    \\ \hline
			\ti & \otimes &  &  
		\end{pmatrix}. 
		\eeq 
	To get a contradiction, suppose the $\otimes$ blocks are  zero. We first redefine $A_j'$ as $A_{j}' := \lambda A_{j}' + A_{j+1}'$ for some large $\lambda >0.$ Then by the definition of positive definiteness (a symmetric matrix $G$ is positive definite if  
	$x^\top  G x > 0$ for all nonzero $x$)
	we find 
	$$
	A_{j}'(\I_{j:(j+1)}) \succ 0.
	$$
	Let $Q$ be a matrix of suitable scaled eigenvectors of $	A_{j}'(\I_{j:(j+1)}), \,$ define 
	
		$$
	T \, := \,  \begin{pmatrix}[c|c|c]
			\bovermat{$\I_{1:(j-1)}$}{\mbox{$\,\,\,\, I  \,\,\,\,$}} 	& \bovermat{$\I_{j:(j+1)}$}{\mbox{$\,\, \pha{\times \times} \,\,$}}	&  \bovermat{$\I_{(j+2):(k+1)}$ }{\mbox{$\,\,\,\,\,\,\pha{\times \times} \,\,\,\,\,\,\,\,\,$}}	\\ \hline 
			 & Q  &    \\ \hline
			 &  & I  
		\end{pmatrix}, 
		$$
		and let 
		$$
		A_i' := T^\top A_i T \; \text{for} \;  i =1, \dots, j, j+2, \dots, k.
		$$
		After this transformation we have that  $A_{j}'(\I_{j:(j+1)}) = I.$ Further, 
		an elementary calculation shows that 
	$(A_1', \dots, A_j', A_{j+2}', \dots, A_k')$ is a length $k-1$ 
	regular facial reduction sequence that satisfies the requirements of Lemma \ref{lemma-reformulation}. However, we assumed that the shortest such sequence has length $k.$ 
	 This contradiction completes the proof.  
\qed  

In Lemma \ref{lemma-derive-poly} we construct a sequence of polynomial inequalities that must be
satisfied  by any $(x_1, \dots, x_k)$ that complete $(\bar{x}_{k+1}, \dots, \bar{x}_m)$ to a strictly feasible solution.
We need some more notation.
Given a strictly feasible solution 
$$ (x_1, \dots, x_k,  \bar{x}_{k+1}, \dots, \bar{x}_m)$$ 
we will write $\delta_j$  for an affine combination of the $``x"$ and $``\bar{x}"$ terms 
with indices larger than $j. \, $ 
In  other words,
\begin{equation} \label{eqn-define-delta} 
\delta_j =\gamma_{j+1} x_{j+1} + \dots + \gamma_k x_k + \gamma_{k+1} \bar{x}_{k+1} + \dots + \gamma_m \bar{x}_m + \gamma_{m+1}, 
\end{equation} 
where the $\gamma_i$ are constants for $i=j+1, \dots, m+1.$  

We will  actually slightly  abuse this notation. We will write $\delta_j$ more than once, but we may mean a different affine combination each time. For example, suppose $k=3, \, $ and $m=4; \, $ 
then we may write $\delta_2 = 2 x_3 + 3 \bar{x}_4 + 5$ on one line, and
$\delta_2 = x_3 - 2 \bar{x}_4 - 3$ on another. 
Given that $\bar{x}_{k+1},  \dots, \bar{x}_m$ are fixed,  $\delta_k$ will always denote a constant.

	\begin{Lemma} \label{lemma-derive-poly} 
		Suppose that $(x_1, \dots, x_k,  \bar{x}_{k+1}, \dots, \bar{x}_m)$ is strictly feasible in \eqref{problem-sdpprime}.  
		Then 
		\begin{equation} \label{eqn-pj} 
		p_j(x_1, \dots, x_k) \, >  \, 0 \,\, \text{for} \,\, j=1, \dots, k-1,
		\end{equation}
		for some $p_j$ polynomials defined as follows: 
			\begin{itemize}
			\item if $t_{j+1} \leq k, \,$ then we choose $p_j$ as 
			\begin{equation} \label{eqn-type1} 
			p_j(x_1, \dots, x_k) \, = \,   (x_j + \delta_j) (x_{t_{j+1}} +  \delta_{t_{j+1}} )  -  ( \beta_{j+1} x_{j+1} + \delta_{j+1})^2,     
			\end{equation}
			\text{where}   $\beta_{j+1}$ is a nonzero constant.
								In this case we call $p_j$ a {\em type 1 polynomial.}

			\item if $t_{j+1} = k+1, \,$ then we choose $p_j$ as 
			\begin{equation} \label{eqn-type2} 
				p_j(x_1, \dots, x_k) \, = \,  	(x_j + \delta_j)     -  ( \beta_{j+1} x_{j+1} + \delta_{j+1})^2,
			\end{equation}
		where $\beta_{j+1}$ is a nonzero constant.
												In this case we call $p_j$ a {\em type 2 polynomial.}
				\end{itemize}

				\qed 

	\end{Lemma}
	Before we prove Lemma  \ref{lemma-derive-poly}, we discuss it. \co{Type 1 polynomials are more difficult to deal with, since they contain 
	mixed  quadratic terms of the form $x_j x_{t_{j+1}}.$ Type 2 polynomials are easier to handle, since they have no such mixed quadratic term. }
First we note that  $p_{k-1}$ will always be 
type 2, since by Lemma    \ref{lemma-tailindex} (with $j=k-1$) we have $t_k = k+1.$ 

In Khachiyan's example \eqref{problem-khachiyan}  all  inequalities come from type 2  polynomials, namely from 
$x_j - x_{j+1}^2 \, $ for $j=1, \dots, k-1.$ 
In contrast, among the inequalities \eqref{eqn-mild-quadratic}  derived from \eqref{problem-mild-sdp} the first 
 two  come from type 1 polynomials and the last one from a type 2 polynomial.

		\pf{of Lemma \ref{lemma-derive-poly}} Fix $j \in \{1, \dots, k-1\}.$ By  the definition of $t_{j+1}, \, $ there is a nonzero element 
		in $A_{j+1}'(\I_j, \I_{t_{j+1}}).$ Let us choose  $\ell_1 \in \I_j$ and $\ell_2 \in \I_{t_{j+1}}$ such that the $(\ell_1, \ell_2)$ element of $A_{j+1}', \, $ which we denote by $(A_{j+1}')_{\ell_1, \ell_2}, $ is nonzero.  

 As stated, suppose that $(x_1, \dots, x_k, \bar{x}_{k+1}, \dots, \bar{x}_m)$ is strictly feasible 
 in \eqref{problem-sdpprime}.  For brevity, define
 \begin{equation}
 S \, := \, \sum_{i=1}^k x_i A_i' +  \sum_{i=k+1}^m \bar{x}_i  A_i'  + B'. 
 \end{equation}

 We distinguish two cases.    
  \begin{enumerate}[wide=18pt]
 	\item [{\bf Case 1:}]  Suppose $t_{j+1} \leq k. $
 	Below we show the matrices that 
 	will be important when we define $p_j:$ 

 
 \begin{small} 
  \beq \label{eqn-S-to-get-pj} 
 \begin{split} 
\underbrace{\begin{pmatrix}[c|c|c|c|c||c]
 		{\mbox{$\,\,\,\,\ti \,\,\,\,$}} 	& \bovermat{$\I_{j}$}{\mbox{$\,\,\, \ti \,\,\,\,$}}	& {\mbox{$\,\,\,\, \ti \,\,\,$}}	& \bovermat{$  \I_{t_{j+1}}$}{\mbox{$\,\,\,\,\ti \,\,\,$}}	 & {\mbox{$\,\,\,\,\ti\,\,\,\,$}} & \bovermat{$\I_{k+1}$}{\mbox{$\,\,\,\,\ti \,\,\,$}}        \\ \hline 
 		\,\,\,\ti \,\, & I & \,\,\,\, & \,\,\,\, & & \\ \hline 
 		\,\, \ti \,\,  & & & & & \\ \hline 
 		\ti & \,\,\,\,  &  &   & &  \\ \hline
 		\ti &  & &   & &  \\ \hline \hline 
 		\ti & \pha{0} &  &  & & 
 	\end{pmatrix}}_{A_{j}'}, & \underbrace{\begin{pmatrix}[c|c|c|c|c||c]
 		{\mbox{$\,\,\,\,\ti \,\,\,\,$}} 	& \bovermat{$\I_{j}$}{\mbox{$\,\,\, \ti \,\,\,\,$}}	& {\mbox{$\,\,\,\, \ti \,\,\,$}}	& \bovermat{$  \I_{t_{j+1}}$}{\mbox{$\,\,\,\,\ti \,\,\,$}}	 & {\mbox{$\,\,\,\,\ti\,\,\,\,$}} & \bovermat{$\I_{k+1}$}{\mbox{$\,\,\,\,\ti \,\,\,$}}        \\ \hline 
 		\,\,\,\ti \,\, & \ti & \ti & \bullet  &  &  \\ \hline 
 		\,\,\,\ti \,\, & \ti & \ti &   &  &  \\ \hline 
 		\,\,\,\ti \,\, & \bullet &  &    &  &  \\ \hline 
 		\ti &  &  & & &  \\ \hline \hline 
 		\ti & &  & & &  \\ 
 	\end{pmatrix}}_{A_{j+1}'}, \,  \\ \\ 
 	  \underbrace{\begin{pmatrix}[c|c|c|c|c||c]
 		{\mbox{$\,\,\,\,\ti \,\,\,\,$}} 	& \bovermat{$\I_{j}$}{\mbox{$\,\,\, \ti \,\,\,\,$}}	& {\mbox{$\,\,\,\, \ti \,\,\,$}}	& \bovermat{$  \I_{t_{j+1}}$}{\mbox{$\,\,\,\,\ti \,\,\,$}}	 & {\mbox{$\,\,\,\,\ti\,\,\,\,$}} & \bovermat{$\I_{k+1}$}{\mbox{$\,\,\,\,\ti \,\,\,$}}        \\ \hline 
 		\,\,\,\ti \,\, & \ti & \ti & \ti & \ti & \ti \\ \hline 
 		\,\,\,\ti \,\, & \ti & \ti & \ti  & \ti & \ti \\ \hline 
 		\,\,\,\ti \,\, & \ti  & \ti & I   &  &  \\ \hline 
 		\ti & \ti & \ti & & &  \\ \hline \hline 
 		\ti & \ti & \ti & & &  \\ 
 	\end{pmatrix}}_{A_{t_{j+1}}'}. 
 \end{split} 
 \eeq
 \end{small} 
As usual, the empty blocks are zero and  the $\ti$ blocks may have arbitrary elements. (More precisely, 
$A_{j+1}'(\I_{j+1})=I,$ but we do not indicate this in equation 
 \eqref{eqn-S-to-get-pj},  since the other entries will suffice to 
 derive the $p_j$ polynomial.) Also,  the $\bullet$ blocks are  nonzero. 
 
 Define  $\beta_{j+1} := (A_{j+1}')_{\ell_1, \ell_2} .$ Let $S'$ be the submatrix of $S$ that contains rows and columns indexed by $\ell_1$ and $\ell_2. \, $ Then 
 $S'$ looks like 
 $$
S' \, = \, \, \bpx x_j + \delta_j                           &  \beta_{j+1} x_{j+1} + \delta_{j+1}  \\
                                                    \beta_{j+1} x_{j+1} + \delta_{j+1}  & x_{t_{j+1}}  + \delta_{t_{j+1}}  \epx.
                                                    $$
                                                    We define $p_j(x_1, \dots, x_k)$ as the determinant of $S', \, $ then 
                                                    $p_j$ is a type 1 polynomial as required in \eqref{eqn-type1}.
                                                      Since $S' \succ 0,$ we see that $p_j(x_1, \dots, x_k)  > 0$ 
                                                      and the proof in this case is complete.

                                                      \item[{\bf Case 2:}] Suppose $t_{j+1} = k+1. \, $ Now $p_j$ will mainly depend on two matrices that we show below: 
                                                      
                                                      \vspace{.1cm} 
                                                        \beq \label{eqn-S-to-pj-2} 
                                                        \begin{split} 
                                                        	\underbrace{\begin{pmatrix}[c|c|c||c]
                                                        			{\mbox{$\,\,\,\,\ti \,\,\,\,$}} 	& \bovermat{$\I_{j}$}{\mbox{$\,\,\, \ti \,\,\,\,$}}	& {\mbox{$\,\,\,\, \ti \,\,\,$}}	&  \bovermat{$  \I_{t_{j+1}}$}{\mbox{$\,\,\,\,\ti \,\,\,$}}	  \\ \hline 
                                                        			\ti & I  &  &  \\ \hline 
                                                        			\ti   &  &  &  \\ \hline \hline 
                                                        			\ti   &  &   &    
                                                        		\end{pmatrix}}_{A_{j}'},  \underbrace{\begin{pmatrix}[c|c|c||c]
                                                        		{\mbox{$\,\,\,\,\ti \,\,\,\,$}} 	& \bovermat{$\I_{j}$}{\mbox{$\,\,\, \ti \,\,\,\,$}}	& {\mbox{$\,\,\,\, \ti \,\,\,$}}	&  \bovermat{$  \I_{t_{j+1}}$}{\mbox{$\,\,\,\,\ti \,\,\,$}}	  \\ \hline 
                                                        		\ti & \ti & \ti & \bullet  \\ \hline 
                                                        		\ti   & \ti & \ti  &   \\ \hline \hline 
                                                        		\ti   & \bullet  &    &    
                                                        	\end{pmatrix}}_{A_{j+1}'}                   \\                     	
                                                        	\end{split} 
                                                        	\eeq  
                                                       Again, the $\bullet$ blocks are nonzero. 
                                                       
                                                      Define  $\lambda := (A_{j+1}')_{\ell_1, \ell_2} \, $  then from the definition of $\ell_1$ 
                                                  	and $\ell_2$ we have $\lambda \neq 0.$ 
                                                  	Let  $\mu :=  S_{\ell_2, \ell_2}, $  then $S \succ 0$ implies $\mu > 0.$ 
                                                  	 Also, since $\ell_2 \in \I_{k+1}, \, $ we 
                                                       see  that $\mu$ depends only on 
                                                       $\bar{x}_{k+1}, \dots, \bar{x}_{m}, $ the $A_i'$ and $B',$ in other words 
                                                       it is a constant.     
                                                       
                                                       We again  let $S'$ be the submatrix of $S$ that contains rows and columns indexed by $\ell_1$ and $\ell_2. $
                                                       Then $S'$ looks like
                                                      $$                                 
                                                       S'  \, = \, \bpx x_j + \delta_j                           & \lambda  x_{j+1} + \delta_{j+1}   \\
                                                      \lambda  x_{j+1} + \delta_{j+1}   & \mu  \epx.
                                                       $$
                                                      Define 
                                                     $$
                                                     p_j(x_1, \dots, x_k) := \dfrac{1}{\mu} \det S'.
                                                     $$
                                                     Since $S' \succ 0, \,$ and  $\mu > 0, \,$ we have $p_j(x_1, \dots, x_k) > 0.$ Thus  
                                                     $$
                                                      p_j(x_1, \dots, x_k)  \, = \, (x_j + \delta_j) - \biggl(  \dfrac{\lambda}{\sqrt{\mu}} x_{j+1} + \dfrac{\delta_{j+1}}{\sqrt{\mu}} \biggr)^2. 
                                                     $$ 
                                                
                                                     Hence  $p_j(x_1, \dots, x_k)$  is a type 2 polynomial 
                                                     in the form required in \eqref{eqn-type2} with 
                                                     $\beta_{j+1} = \lambda/\sqrt{\mu}.$ (Since $\mu$ is a constant, by our definition of $\delta_{j+1}$ in \eqref{eqn-define-delta} 
                                                     we have that 
                                                       $\delta_{j+1}/\sqrt{\mu}$ is still $\delta_{j+1}.$) 
                      The proof in this case is now complete.

                                                     \end{enumerate}                                                        
                                                    
                                                       \qed

                                                       \begin{Lemma} \label{lemma-xj>=dj+1} 
                                                       	Suppose that $(x_1, \dots, x_k, \bar{x}_{k+1}, \dots, \bar{x}_m)$ 
                                                       	is strictly feasible in \eqref{problem-sdpprime} and $x_k$ is sufficiently large.
                                                       	                                                       	Then 
                                                       		\begin{equation} \label{eqn-xjgeqdj+1} 
                                                       			x_j \geq d_{j+1} x_{j+1}^{\alpha_{j+1}} \,\, \text{for \,} j=1, \dots, k-1,
                                                       		\end{equation}  
                                                       		where  the $d_{j+1}$ are positive constants and the 
                                                       		$\alpha_{j+1}$ satisfy the recursion 
                                                   		\begin{equation} \label{eqn-alpha-j-recursion}  
                                                       		\alpha_{j+1} \, = \, \left\{   \begin{array}{rll} 
                                                       		2 - \dfrac{1}{\alpha_{j+2} \dots \alpha_{t_{j+1}}} & \text{if} &  t_{j+1} \leq k \\ 
                                                       		 2 & \text{if} &  t_{j+1} = k+1 
                                                       		\end{array}  \right. 
                                                       		\end{equation} 
                                                       		for $j=1, \dots, k-1.$ 
                                                       			                                              		                                                    		
                                                       	\qed
                                                       \end{Lemma}
                                                       
                                                       Before we prove  Lemma \ref{lemma-xj>=dj+1}, we discuss it.
                                                       We have $t_k = k+1$ (by Lemma  \ref{lemma-tailindex}) hence  Lemma \ref{lemma-xj>=dj+1} 	implies $\alpha_k=2. \, $ 
                                                       Hence,  by induction the recursion  \eqref{eqn-alpha-j-recursion} implies that 
                                                       $\alpha_j \in (1,2]$ holds for all $j$ (naturally, we compute $\alpha_k, \alpha_{k-1}, \dots, \alpha_2$ in this order).  
                                                       Thus, if  $x_k$ is large enough, then $x_j > 0$ for  $j=1, \dots, k.$

                                                       It is also interesting that  formula \eqref{eqn-alpha-j-recursion}  
                                                       is reminiscent of a continued fractions formula. 
                                                       
                                                       To illustrate Lemma \ref{lemma-xj>=dj+1} we show how from 
                                                       \eqref{problem-mild-sdp} we can deduce the inequalities 
                                                         \eqref{eqn-mild-quadrati-corollary} much more quickly than we did before. Recall that in this example we have $k=4.$
                                              We compute 
                                                         the $\alpha_{j+1}$ exponents by the recursion  \eqref{eqn-alpha-j-recursion} 
                                                         as 
                                                        \begin{equation} \label{eqn-alpha-j-rec-example} 
                                                        \begin{array}{rclr} 
                                                        \alpha_4 & = & 2                                         & (\text{since} \; t_4 = 5) \\
                                                        \alpha_3 & = & 2 - 1/\alpha_4 \, = \, 3/2 & (\text{since} \; t_3 = 4) \\
                                                        \alpha_2 & = & 2 - 1/\alpha_3  \, = \, 4/3 & \,\,\,\,\, (\text{since} \; t_2 = 3).
                                                        \end{array}
                                                        \end{equation}
                                                      
                                                      Next we  sketch the proof of Lemma \ref{lemma-xj>=dj+1}. We start with the inequalities 
                                                      $p_j(x_1, \dots, x_k) > 0$ derived in Lemma \ref{lemma-derive-poly}; 
                                                      these are satisfied by all strictly feasible solutions of \eqref{problem-sdpprime}. Note that the $p_j$ polynomials defined in  \eqref{eqn-type1} and 
                                                       \eqref{eqn-type2}  are quite messy. However, if $p_j$ is a type 1 polynomial (defined in \eqref{eqn-type1}), then we 
                                                       deduce a 
                                                    cleaned up inequality 
                                                       $$
                                                       x_j x_{t_{j+1}} > \const x_{j+1}^2,
                                                       $$
                                                       assuming $x_k$ is large enough. 
                                                       Similarly, if $p_j$ is a type 2 polynomial (defined in \eqref{eqn-type2}), 
                                                       then we derive a similarly cleaned up inequality
                                                        $$
                                                       x_j  >  \const x_{j+1}^2,
                                                       $$
                                                       assuming $x_k$ is large enough. 
                                                       Then from 
                                                       the cleaned up inequalities we derive the required inequalities \eqref{eqn-xjgeqdj+1} and the recursion 
                                                        \eqref{eqn-alpha-j-recursion}.

                                      Next, since the proof of Lemma \ref{lemma-xj>=dj+1} is somewhat technical, we illustrate the cleaning up  of the inequalities with  an example.

                                                                               \begin{figure} 
                                                                               \begin{center}
                                                                               \co{\includegraphics[scale=0.5]{pics/rainbow/K_v1.png}	\includegraphics[scale=0.5]{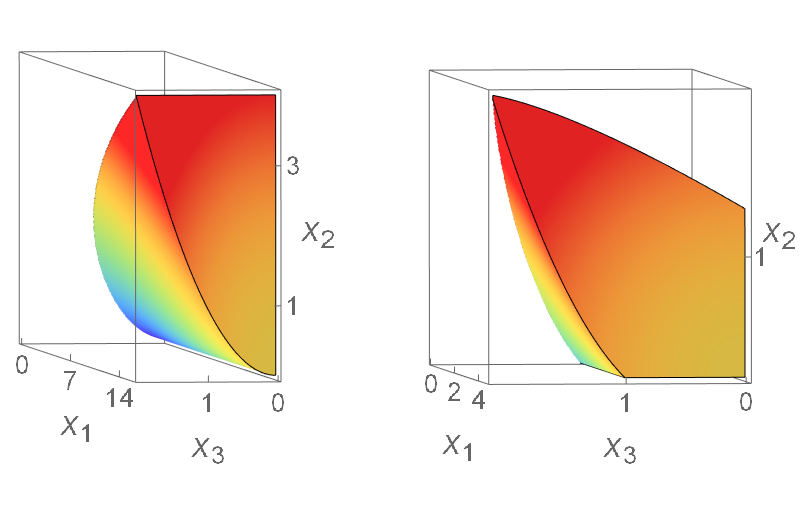}
                                                                            }
                                                                               \co{
                                                                               	
                                                                               	\includegraphics[scale=0.5]{pics/rainbow/K_messy_v1.png}
                                                                               }
                                                                               	
                                                                               	\includegraphics[scale=0.5]{pics/rainbow/K_v1.png}
                                                                               	\includegraphics[scale=0.5]{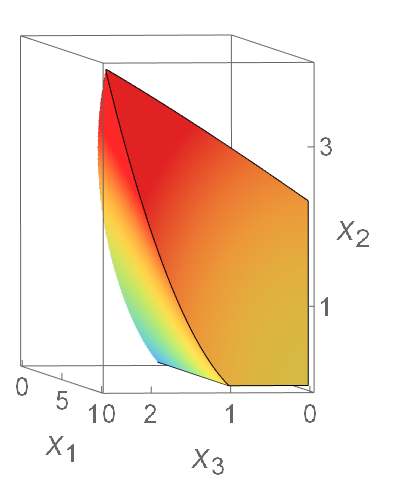}

                                                                                \end{center} 
                                                                               \caption{Feasible sets of \eqref{problem-khachiyan} (on left) and of the inequalities derived from the perturbed Khachiyan SDP \eqref{problem-khachiyan-messy}
                                                                               	 (on the right)} 
                                                                               \label{figure-khach-vs-messy} 
                                                                                  \end{figure}

                                                  \begin{Example} (Perturbed Khachiyan)  \label{example-khach-messy} 
                                                  	As a warmup, first let us consider the SDP  
                                                  		\renewcommand{\arraystretch}{1.4} 
                                                  	\begin{equation} \label{problem-khachiyan-3} 
                                                  		\begin{pmatrix}[cccc]  x_1  &                       &      & x_2  \\ 
                                                  			& x_2       &      & x_3 \\ 
                                                  			&                        &  x_3 &     \\  
                                                  			x_2             &  x_3                 &         & 1 \end{pmatrix} \succeq 0, 
                                                  	\end{equation}
                                                  which is just a smaller version of   \eqref{problem-khach-sdp}. Thus the feasible solutions of \eqref{problem-khachiyan-3} 
                                                  satisfy the inequalities 
                                                  \begin{equation} \label{eqn-khach-ineq-3} 
                                                  	x_1 \geq x_2^2, \, x_2 \geq x_3^2. 
                                                  \end{equation}
                                                  	Next, let us consider a ``perturbed" version of \eqref{problem-khachiyan-3}
                                                  	\renewcommand{\arraystretch}{1.4} 
                                                  	\begin{equation} \label{problem-khachiyan-messy} 
                                                  	\begin{pmatrix}[c|c|c|c]  x_1 -2x_2 &                       &      & x_2 - x_3 \\ \hline 
                                                  	& x_2 + x_3       &      & x_3 \\ \hline 
                                                  	&                        &  x_3 &     \\ \hline 
                                                  	x_2 - x_3            &  x_3                 &         & 1 \end{pmatrix} \succeq 0,
                                                  	\end{equation}
                                                  which we obtain from \eqref{problem-khachiyan-3} 
                                                  by  replacing 
                                                  $x_1$ by $x_1 -2x_2$ and $x_2$ by $x_2 \pm x_3.$ 
                                                  
                                                	Suppose   $(x_1, x_2, x_3)$ is feasible in \eqref{problem-khachiyan-messy}
                                                	Then  the quadratic inequalities \eqref{eqn-khach-ineq-3} may not hold  \footnote{
                                                	For example, $x=(5,2,2)$ is feasible 
                                                		in  \eqref{problem-khachiyan-messy}, but not in \eqref{eqn-khach-ineq-3}.}. 
                                                  	However, let us also assume  $x_3 \geq 10.$ We then claim that the 
                                                  	following slightly weaker inequalities do hold: 
                                                  \begin{eqnarray}  \label{eqn-x1-perturbed}  
                                                  	x_1 & \geq & \frac{1}{2} x_2^2 \\
                                                  	\label{eqn-x2-perturbed} 
                                                  	x_2 & \geq & \frac{1}{2} x_3^2.
                                                  \end{eqnarray}
                                              To prove that, from    the                       
                                       principal minors of \eqref{eqn-khach-ineq-3} we first deduce the inequalities 
                                                  	\begin{eqnarray} \label{eqn-x1}
                                                  	x_1 -2x_2 & \geq & (x_2 - x_3)^2 \\ \label{eqn-x2}
                                                  	x_2 + x_3 & \geq & x_3^2.
                                                  	\end{eqnarray}
                                                  Then  \eqref{eqn-x2-perturbed} follows from \eqref{eqn-x2} and $x_3 \geq 10$ directly. Hence $x_2 \geq 50$ also holds. 
                                                  
                                                  To prove \eqref{eqn-x1-perturbed}, we lower bound the right hand side in \eqref{eqn-x1} as
                                                  \begin{equation} 
                                                  	\begin{array}{rcl} 
                                                  		(x_2 - x_3)^2 & \geq & (x_2 - \sqrt{2 x_2})^2 \\
                                                  		& \geq & \dfrac{1}{2} x_2^2,
                                                  		\end{array} 
                                                  	\end{equation} 
                                                  where the first inequality is from  $x_2 \geq x_3$ and  \eqref{eqn-x2-perturbed}. The second inequality follows, since $x_2 \geq 50.$ 
                                                  Using this lower bound in \eqref{eqn-x1} and $x_2 \geq 0, \, $ the desired 
                                                  inequality  \eqref{eqn-x1-perturbed} follows.
                                               \end{Example} 
                                                       	We show the  set described by the inequalities 
                                                       \eqref{eqn-khach-ineq-3} which appear in \eqref{problem-khachiyan}, 
                                                       	and the feasible set  
                                                       	described by the inequalities \eqref{eqn-x1}-\eqref{eqn-x2}  on Figure 
                                                       	\ref{figure-khach-vs-messy}. We normalized both sets by suitable bounds on $x_3.$ 
                                                       	 Note that $x_3$ increases from right to left for better visibility.

                                                       \pf{of Lemma   \ref{lemma-xj>=dj+1}} We use an argument analogous to the one in Example  \ref{example-khach-messy}. 
                                                       We use induction and show how to suppress  the $``\delta"$ terms in the type 1 and type 2 polynomials at the cost of  making $x_k$ large and choosing suitable $d_j$ constants. 
                                                       
                                                       
                                                       Suppose that $(x_1, \dots, x_k, \bar{x}_{k+1}, \dots, \bar{x}_m)$ is strictly feasible 
                                                       in \eqref{problem-sdpprime}. Then by Lemma \ref{lemma-derive-poly} the inequalities $p_j(x_1, \dots, x_k) > 0$ hold for $j=1, \dots, k-1.$ 
                                                       
                                                     We first establish the base case for the induction. 
                                                  	From the remark after the statement of Lemma \ref{lemma-derive-poly} we recall that 
                                                  	$p_{k-1}$ is a type 2 polynomial. Hence 
                                                        \begin{equation} \label{kramer-pre} 
                                                    \begin{array}{rcl}
                                                     p_{k-1}(x_1, \dots, x_k) & = &  (x_{k-1} + \delta_{k-1})   - ( \beta_k x_k + \delta_k)^2 \\
                                                    	    & = &  (x_{k-1} + \gamma_k x_k + \delta_{k})   - ( \beta_k x_k + \delta_k)^2 \,  > \, 0,
                                                    	\end{array} 
                                                     \end{equation}
                                                   where the first equality is from the definition of $p_{k-1} \, $  (see \eqref{eqn-type2} with $j=k-1. \, $). Here 
                                                   $\beta_k \neq 0$ is a constant. The 
                                                   second equality is from the definition of $\delta_{k-1}$ (see \eqref{eqn-define-delta}), 
                                                   where $\gamma_k$ is a constant which may be zero. 
                                                   
                                                   Since $\delta_k$ is a constant, and $\beta_k \neq 0, \,$ 
                                                   from \eqref{kramer-pre} we deduce that 
                                                    $x_{k-1} \geq d_k x_k^2$ if $x_k$ is sufficiently large, where  $d_k$ is a suitable positive constant. So the proof of the base case is complete.

                                                       For the inductive step we will adapt the $O, \Theta$ and $o$ notation 
                                                       from theoretical computer science. 
                                                       Given functions $f, g: \rad{k}  \rightarrow \rad{}_+$ we say that 
                                                       \begin{enumerate}
                                                       	\item $f=O(g)$ (in words, $f$ is big-Oh of $g$) if there are positive constants $C_1$ and $C_2$ such that 
                                                       	if  $(x_1, \dots, x_k, \bar{x}_{k+1}, \dots, \bar{x}_{m})$ is strictly feasible in \eqref{problem-sdpprime} and $x_k \geq  C_1, \,$ then 
                                                       	$$
                                                       	f(x_1, \dots, x_k) \leq C_2 g(x_1, \dots, x_k). 
                                                       	$$
                                                       	\item $f = \Theta(g)$ (in words, $f$ is big-Theta of $g$) 
                                                       	if $f = O(g)$ and $g=O(f).$ 
                                                       	\item $f=o(g)$ (in words, $f$ is little-oh of $g$) 
                                                       	 if for all  $\epsilon > 0$  there is $\delta > 0$  such that 
                                                       	if 	$(x_1, \dots, x_k, \bar{x}_{k+1}, \dots, \bar{x}_{m})$ is strictly feasible in \eqref{problem-sdpprime} and $x_k \geq  \delta$ then 
                                                       		$$
                                                       	f(x_1, \dots, x_k) \, \leq   \, \epsilon  g(x_1, \dots, x_k). 
                                                       	$$
                                                       \end{enumerate}
                                                   
                                          The usual calculus of $O, \Theta$ and $o$ carries over verbatim. We spell out one calculus rule 
                                                       that we will use repeatedly: 
                                                    \begin{equation} \label{eqn-absorb} 
                                                   	\begin{array}{rcl}
                                                   		f=o(h), g = \Theta(h) & \Rightarrow & f+g= \Theta(h).
                                                   	\end{array} 
                                                   \end{equation}
In the implication \eqref{eqn-absorb} we say informally  that we {\em absorb}  the $o(h)$ term into the $\Theta(h)$ term. 

For brevity, in the rest of this proof we will say that $(x_1, \dots, x_k)$ is {\em good}, if $(x_1, \dots, x_k, \bar{x}_{k+1}, \dots, \bar{x}_{m})$ is 
strictly feasible in \eqref{problem-sdpprime} and $x_k$ is sufficiently large.

                                                                    Suppose next that $1 <  j+1 \leq k-1$ and we have proved the following:
                                                                    for all good $(x_1, \dots, x_k)$ 
                                                                                                                                    the inequalities 
                                                                     \begin{equation} \label{eqn-xj+1dj+1} 
                                                                     \ba{rcl} 
                                                                     x_{j+1} & \geq & d_{j+2} x_{j+2}^{\alpha_{j+2}} \\
                                                                     x_{j+2} & \geq & d_{j+3} x_{j+3}^{\alpha_{j+3}} \\
                                                                                   & \vdots & \\
                                                                                   x_{k-1} & \geq & d_{k} x_{k}^{\alpha_{k}}
                                                                                   	\end{array}
                                                                           \end{equation}
                                                                  hold,  where $d_{j+1}, \dots, d_k$ are 
                                                                   positive constants, and  $\alpha_{j+2}, \dots, \alpha_k$ are positive constants  
                                                                  derived from the recursion 
                                                                  \eqref{eqn-alpha-j-recursion}; further,  
																	$\alpha_k = 2.$

                                                         	We will next show that for all good  $(x_1, \dots, x_k)$
                                                                                                                            the inequality 
                                                                           \begin{equation} \label{xj-toget} 
                                                                            x_{j} \,  \geq \,  d_{j+1} x_{j+1}^{\alpha_{j+1}} 
                                                                           \end{equation}
                                                                           holds, where   $d_{j+1} \, $ is another positive constant and 
                                                                           $\alpha_{j+1}$ is computed by  the recursion \eqref{eqn-alpha-j-recursion}.

                                                      For that, first note that 
                                                      the recursion  \eqref{eqn-alpha-j-recursion} implies by straightforward 
                                                       induction that 
                                                       $\alpha_{k}, \alpha_{k-1},$ $ \dots, \alpha_{j+2}$ are in the interval $(1,2].$
                                                        So by the inequalities 
                                                     \eqref{eqn-xj+1dj+1} we have 
                                                          	\begin{equation} \label{eqn-xj-o} 
                                                       	x_{s}  =  \, o(x_{\ell}) \,\,\, \text{when} \,\,  s > \ell \geq j+1.   
                                                       	\end{equation} 
                                                       	                              
                                                       We distinguish two cases.                  	                                              	   
                                                       	   \begin{enumerate}[wide=18pt]
                                                       	   	\item[{\bf Case 1:}]
                                                       	         	 	  First suppose that $t_{j+1} \leq k, \, $ in other words,  
                                                       	         	    the quadratic polynomial $p_j$ is type 1
                                                       	         	     (see \eqref{eqn-type1}).  
                                                       	         	      Then we  claim that for all good $(x_1, \dots, x_k)$ the following hold: 
                                                       	         	      	    \begin{equation} \label{eqn-jerry}   
                                                       	         	    	\begin{split} 
                                                       	         	    		0 & < p_j(x_1, \dots, x_k)  \\
                                                       	         	    		    & = (x_j + \delta_j ) (x_{t_{j+1}} + \delta_{t_{j+1}}) - (\beta_{j+1} x_{j+1}  + \delta_{j+1})^2  \\
                                                       	         	    		   	& =  (x_j + \gamma_{j+1} x_{j+1} + \delta_{j+1} ) (x_{t_{j+1}} + \delta_{t_{j+1}}) - (\beta_{j+1} x_{j+1}  + \delta_{j+1})^2 \\
                                                       	         	    		   	& \leq   (x_j + |\gamma_{j+1}|  x_{j+1} + |\delta_{j+1}| ) (x_{t_{j+1}} + |\delta_{t_{j+1}}|) - (\beta_{j+1} x_{j+1}  + \delta_{j+1})^2 \\
                                                       	         	    		& =  (x_j + |\gamma_{j+1}|  x_{j+1} + o(x_{j+1}))  (x_{t_{j+1}} + o(x_{t_{j+1}}))   - (\beta_{j+1} x_{j+1}  + o(x_{j+1}))^2.
                                                       	         	    	\end{split}
                                                       	         	    \end{equation}
                                                                	    Here $\beta_{j+1}$ is a nonzero constant and  $\gamma_{j+1}$ is a constant which may be zero. 
                                                                	    	
                                                                	    	Indeed, in \eqref{eqn-jerry} the first inequality is by Lemma \ref{lemma-derive-poly} and 
                                                                	    	the first equality is from the definition of $p_j$ in \eqref{eqn-type1}.
                                                                	   The second equality follows 
                                                                	    from the definition of $\delta_{j}$ in \eqref{eqn-define-delta}.   
                                                                	    The second inequality follows, since $j+1 > 1, $ hence by Lemma \ref{lemma-tailindex} we have 
                                                                	    $t_{j+1} > j+1, \, $ so for good $(x_1, \dots, x_k)$ 
                                                                	    by the inequalities \eqref{eqn-xj+1dj+1}  we have 
                                                                	    $x_{t_{j+1}} > 0$ and  $x_{t_{j+1}} + \delta_{t_{j+1}} > 0.$ 
                                                                	    	 The third equality follows, since 
                                                                	    by \eqref{eqn-xj-o} the term $\delta_{j+1}$ is a linear combination of $o(x_{j+1})$ terms; and 
                                                                	    by $t_{j+1} > j+1$  and by  \eqref{eqn-xj-o} the term 
                                                                	    $\delta_{t_{j+1}}$ is a linear combination of $o(x_{t_{j+1}})$ terms.

                                                                	    We next claim that the last expression in \eqref{eqn-jerry} is upper bounded by 
                                                                	    \begin{equation} \label{eqn-jerry-2} 
                                                                	    	(x_j + \Theta (x_{j+1}) ) \Theta(x_{t_{j+1}}) - \Theta(x_{j+1})^2. 
                                                                	    \end{equation}
                                                              For that, first assume $\gamma_{j+1} \neq 0.$ Then 
                                                            by the rule \eqref{eqn-absorb} we absorb the $o(x_{j+1})$ terms into the terms with $x_{j+1};$ and the $o(x_{t_{j+1}})$ 
                                                            term into the term with $x_{t_{j+1}}.$  Next, assume $\gamma_{j+1} = 0.$ Then we use the bound 
                                                            $0 =  \gamma_{j+1} x_{j+1}  < \Theta(x_{j+1}), \, $ 
                                                            and for the rest of the estimate we still use the absorbing rule
                                                            \eqref{eqn-absorb}. 
                                                            
                                                       	    	We then continue \eqref{eqn-jerry} by using the upper bound \eqref{eqn-jerry-2}:  
                                                       	    	\begin{equation} \label{eqn-jerry-3} 
                                                       	    	\begin{split} 
                                                       	    	0 & < (x_j + \Theta (x_{j+1}) ) \Theta(x_{t_{j+1}}) - \Theta(x_{j+1})^2  \\
                                                       	    	   & = x_j  \Theta(x_{t_{j+1}}) + \Theta (x_{j+1} x_{t_{j+1}})  - \Theta(x_{j+1})^2  \\
                                                       	    	   & \leq  x_j \Theta( x_{t_{j+1}}) + o(x_{j+1}^2)  - \Theta(x_{j+1})^2 \\
                                                       	    	   & \leq   x_j \Theta( x_{t_{j+1}})  - \Theta(x_{j+1})^2, 
                                                       	    	\end{split} 
                                                       	    	\end{equation}
                                                       	    	where the second inequality follows, since $t_{j+1} > j+1$ hence by 
                                                       	    	\eqref{eqn-xj-o} we have $x_{t_{j+1}}=o(x_{j+1}).$ 
                                                       	    	For the last inequality we absorb the 
                                                       	    	$o(x_{j+1}^2)$ term into the $\Theta(x_{j+1}^2)$ term by the rule \eqref{eqn-absorb}. 

                                                       	       	Next we use  $t_{j+1} > j+1$ and combine the inequalities \eqref{eqn-xj+1dj+1} to 
                                                       	    	 learn that 
                                                       	    	\begin{equation} \nonumber 
                                                       	    	x_{t_{j+1}}^\alpha = O(x_{j+1}) \, 
                                                       	    	\end{equation}
                                                       	    	where $\alpha = \alpha_{j+2} \alpha_{j+3} \dots \alpha_{t_{j+1}}.$ 
                                                       	    	Hence $x_{t_{j+1}}  = O(x_{j+1}^{1/\alpha}).$ 
                                                       	    	
                                                       	    	We next plug this last estimate into the last inequality in  \eqref{eqn-jerry-3} and  deduce
                                                       	    	$$
                                                       	    	0 <  x_j 	\Theta(x_{j+1}^{1/\alpha}) - \Theta (x_{j+1}^2).
                                                       	    	$$
                                                       	    	We finally divide this last inequality by $x_{j+1}^{1/\alpha}$ and a constant, and deduce that 
                                                       	    	$$
                                                       	    	x_j \geq d_{j+1} x_{j+1}^{2 - 1/\alpha}
                                                       	    	$$
                                                       	    for a suitable positive $d_{j+1}$ constant, 	if $x_k$ is large enough. Hence we can set 
                                                       	    $\alpha_{j+1} := 2 - 1/\alpha$ and the proof in this case is complete.

                                                       	    		\item[{\bf Case 2}] Suppose that $t_{j+1}=k+1, \, $ in other words, the quadratic polynomial $p_j$ is type $2.$ Then by 
                                                       	    		Lemma \ref{lemma-derive-poly} we get 
                                                       	    		$$
                                                       	    		\begin{array}{rcl} 
                                                       	    		0 < p_j(x_1, \dots, x_k) & = & (x_j + \delta_j)   - ( \beta_{j+1} x_{j+1} + \delta_{j+1})^2 \\
                                                       	    		& = & (x_j + \gamma_{j+1} x_{j+1} + \delta_{j+1})    - ( \beta_{j+1} x_{j+1} + \delta_{j+1})^2 
                                                       	    		\end{array}
                                                       	    		$$
                                                       	    		for some  $\beta_{j+1} \neq 0$ and $\gamma_{j+1}$ constants, where $\gamma_{j+1}$ may be zero.       
                                                       	    		Here the second equality is from the definition of $\delta_j$  in \eqref{eqn-define-delta}.   
                                                       	    		By \eqref{eqn-xj-o}  we have $ \delta_{j+1} = o(x_{j+1}),$  
                                                       	    		hence 
                                                       	    		$$
                                                       	    		x_j  \geq d_{j+1} x_{j+1}^2
                                                       	    		$$
                                                       	    		for a suitable positive constant $d_{j+1}$ if $x_k$ is large enough. 
                                                       	    		So we can set $\alpha_{j+1}=2, $ and the proof is complete. 
                                                       	    			
                                                       	    			 \end{enumerate}
                                                       	    		
                                                       	    		\qed

                                                       	   As a prelude to Lemma \ref{lemma-alphaj-monotone}, in Figure 
                                                       	   \ref{figure-shift} we show three SDPs (for brevity we left out the $\succeq$ symbols).
                                                       	   The first is \eqref{problem-mild-sdp}. The second and third arise from it 
                                                       	   by shifting $x_2$ in the upper offdiagonal position to the right.
                                                       	   Underneath we show the vector of the $\alpha = (\alpha_2, \alpha_3, \alpha_4)$ exponents in the inequalities derived by the recursion \eqref{eqn-alpha-j-recursion}.

                                                       	   We see that $\alpha_2$  
                                                       	   increases from left to right and Lemma \ref{lemma-alphaj-monotone} presents a general result  of this kind.

                                                       	   \begin{figure} 
                                                       	   	\begin{center}
                                                       	   		$\underbrace{\begin{pmatrix}
                                                       	   			\mrbmath{x_1} 	 &           &  \mrbmath{x_2}          &     &   \\
                                                       	   					& \mbbmath{x_2}   &           &  \mbbmath{x_3}   &   \\
                                                       	   					\mrbmath{x_2} &           & \mgbmath{x_3}     &     & \mgbmath{x_4}  \\
                                                       	   					& \mbbmath{x_3}         &            & \mblbmath{x_4}   &         \\
                                                       	   					&   & \mgbmath{x_4}          &  & 	\mblbmath{1}       
                                                       	   				\end{pmatrix}}_{\mblbmath{\alpha=}\mblbmath{(}\mrbmath{4/3}, \,\,\, \mbbmath{3/2}, \,\,\, \mgbmath{2}\mblbmath{)}} 
                                                       	   			$
                                                       	   			 \mblbmath{\rightarrow} 
                                                       	   			$
                                                       	   				\underbrace{\begin{pmatrix}
                                                       	   					\mrbmath{x_1} 	 &           &            &  \mrbmath{x_2}   &   \\
                                                       	   						& \mbbmath{x_2}   &           &  \mbbmath{x_3}   &   \\
                                                       	   						&           & \mgbmath{x_3}     &     & \mgbmath{x_4}  \\
                                                       	   						\mrbmath{x_2} & \mbbmath{x_3}         &            & \mblbmath{x_4}   &         \\
                                                       	   						&   & \mgbmath{x_4}          &  & 	\mblbmath{1}       
                                                       	   					\end{pmatrix}}_{\mblbmath{\alpha=}\mblbmath{(}\mrbmath{5/3}, \,\,\, \mbbmath{3/2}, \,\,\, \mgbmath{2}\mblbmath{)}} 
                                                       	   				$ 
                                                       	   				\mblbmath{\rightarrow}  
                                                       	   				$ 
                                                       	   				\underbrace{\begin{pmatrix}
                                                       	   						\mrbmath{x_1} 	 &           &            &     &  	\mrbmath{x_2}   \\
                                                       	   						& \mbbmath{x_2}   &           &  \mbbmath{x_3}   &   \\
                                                       	   						&           & \mgbmath{x_3}     &     & \mgbmath{x_4}  \\
                                                       	   						& \mbbmath{x_3}         &            & \mblbmath{x_4}   &         \\
                                                       	   						\mrbmath{x_2}  &   & \mgbmath{x_4}          &  & 	\mblbmath{1}       
                                                       	   					\end{pmatrix}}_{\mblbmath{\alpha=}\mblbmath{(}\mrbmath{2}, \,\,\, \mbbmath{3/2}, \,\,\, \mgbmath{2}\mblbmath{)}} 
                                                       	   				$
                                                       	   			\end{center}
                                                       	   		\caption{Shifting $x_2$ to the right increases $\alpha_2$} 
                                                       	   			\label{figure-shift} 
                                                       	   		\end{figure}

                                                       	   	    \begin{Lemma} \label{lemma-alphaj-monotone}  The $\alpha_j$ exponents in 
                                                       	   	    		\eqref{eqn-alphaj-bounds} are strictly  increasing  functions of the 
                                                       	   	    		$t_{j+1}$ tail-indices defined in Definition \ref{definition-tailindex}. 
                                                       	   	    		
                                                       	   	    	Precisely, suppose we derived the inequalities 
                                                       	   	    		\begin{equation} \label{eqn-xj-dj+1-repeat} 
                                                       	   	    		x_\ell \geq d_{\ell+1} x_{\ell+1}^{\alpha_{\ell+1}} \,\, \text{for} \,\, \ell=1, \dots, k-1
                                                       	   	    		\end{equation} 
                                                       	   	    	 from \eqref{problem-sdpprime} using the recursion \eqref{eqn-alpha-j-recursion}.
                                                       	   	    	 	Here $d_{\ell+1}$ is a positive constant for all $\ell.$ 
                                                       	   	    		
                                                       	   	    		Let $j$ be an index in  $\{1, \dots, k-1\} \, $ such that $t_{j+1} \leq k. \,$ 
                                                       	   	    		Suppose we increase $t_{j+1}$ by $1$ (by changing $A_{j+1}'$ in \eqref{problem-sdpprime}), 
                                                       	   	    		then  
                                                       	   	    				derive the  inequalities 
                                                       	   	    	    		\begin{equation} \label{eqn-xj-fj+1} 
                                                       	   	    	    		x_\ell \geq f_{\ell+1} x_{\ell+1}^{\omega_{\ell+1}} \;\, \text{for} \;\, \ell=1, \dots, k-1,
                                                       	   	    	    		\end{equation}
                                                       	   	    	    		  using the recursion \eqref{eqn-alpha-j-recursion}. 
                                                       	   	    	    		Here $f_{\ell+1}$ is a positive constant for all $\ell.$ 
                                                       	   	    	    		
                                                       	   	    	    		Then 
                                                       	   	    	    		\begin{equation}
                                                       	   	    	    		\omega_{\ell+1}  \, \left\{   \begin{array}{rl}  
                                                       	   	    	    		=  \alpha_{\ell+1}  & \text{if} \,\, \ell > j  \\
                                                       	   	    	    		>   \alpha_{\ell+1} & \text{if} \,\, \ell =  j  \\
                                                       	   	    	    		\geq   \alpha_{\ell+1}  & \text{if} \,\, \ell <  j. 
                                                       	   	    	    		\end{array}                                      	   	    	    		                                                       	   	    	    		\right.
                                                       	   	    	    		\end{equation}
                                                       	   	    	    		
                                                       	   	    \end{Lemma}

                                                     \pf{}   Let us make all the assumptions  and note that $t_{j+1} \leq k$  implies that polynomial $p_j$ is type 1.
                                                     
                                                      We first prove $\omega_{\ell+1} = \alpha_{\ell+1}$ for all  $\ell > j.$     For that,                                          	
                                                  	we  observe two  facts. First, when we define the polynomials in Lemma \ref{lemma-derive-poly},  
                                                  	the only polynomial that refers to  $t_{j+1}$ is $p_j.$ 
                                                   Second,  by the proof of Lemma \ref{lemma-xj>=dj+1}, $p_j$ is only used to derive 
                                                      the inequality \eqref{xj-toget}, and hence to determine the value of $\alpha_{j+1}.$ Of course, $\alpha_{j+1}$ affects 
                                                      $\alpha_j, \alpha_{j-1}, \dots, \alpha_2$ via the recursion \eqref{eqn-alpha-j-recursion}.
                                                      However, $\alpha_{j+1}$ 
                                                       does not affect             $\alpha_{\ell+1}$ for $\ell > j.$           
                                                                                                         From these two  facts our claim follows.

                                                     We next prove $\omega_{j+1} > \alpha_{j+1}.$ 
                                                     For brevity, let $s := t_{j+1}$ and 
                                                     $\alpha := \alpha_{j+2} \cdot \dots \cdot \alpha_{s}.$ 
                                                     We first observe that since $s \leq k, \,$ the recursion formula \eqref{eqn-alpha-j-recursion}  shows 
                                                         \begin{eqnarray} \nonumber 
                                                         \alpha_{j+1} & = & 2 - \dfrac{1}{\alpha}. 
                                                         \end{eqnarray}
                                                     We next examine how we 
                                                     compute $\omega_{j+1}$ by formula \eqref{eqn-alpha-j-recursion}.
                                                     We distinguish two cases. If $s < k, \, $ then  in  \eqref{eqn-alpha-j-recursion} 
                                                      we use the top equation, so we get 
                                                     $\omega_{j+1} = 2 - 1/(\alpha \cdot \alpha_{s+1}).$ Since  $\alpha_{s+1} > 1, \,$ we deduce 
                                                     $\omega_{j+1} > \alpha_{j+1}, \, $ as wanted. If  
                                                     $s = k, \, $ then $s+1 = k+1, \,$ so in  \eqref{eqn-alpha-j-recursion}  we use the bottom equation.                                                     
                                                     Thus $2 = \omega_{j+1}$ and 
                                                 $2 > \alpha_{j+1}, \,$  so  
                                                     $\omega_{j+1} > \alpha_{j+1}$ again follows. 
                                                    
                                               Finally,  we prove  $\omega_{\ell+1} \geq \alpha_{\ell+1}$ for $\ell >   j. $ Since we already proved 
                                               this relation for all $\ell \leq  j, \, $ our claim follows 
                                                   by induction from the recursion formula \eqref{eqn-alpha-j-recursion}. 
                                                     \qed

                                                Recall from Lemma \ref{lemma-tailindex} that the tail-index $t_{j+1}$ is at least $j+2$ for all $j.$ In our final lemma 
                                                we examine the case when $t_{j+1}$ equals $j+2$ for all $j$ and we derive a closed form solution for the $\alpha_{j+1}$ exponents. 

                                                    \begin{Lemma} \label{lemma-alphaj-smallest} 
                                                    	Suppose that $t_{j+1} = j+2$ for $j=1, \dots, k-1.$ 
                                                    Then the recursion formula \eqref{eqn-alpha-j-recursion}  yields 
                                                     \begin{equation} \label{eqn-alpha-j+1-k} 
                                                    \alpha_{j+1} = 1 + \dfrac{1}{k-j} \; \text{for} \; j=1, \dots, k-1. 
                                                    \end{equation}
                                                    \end{Lemma}
                                                    \pf{}  We use induction. First suppose $j=k-1.$ Since 
                                                    $p_{k-1}$ is of type 2, we see $\alpha_{j+1} = \alpha_{k}=2,$ as wanted. 
                                                                                                        Next assume that $1 \leq j <  k -1$ and 
                                                   $$
                                                   \alpha_{j+2} = 1 + \dfrac{1}{k-j-1}.
                                                    $$
                                                    \co{$
                                                    \alpha_{j+2} = 1 + 1/(k-j+1). 
                                                    $}
                                                   By the recursion \eqref{eqn-alpha-j-recursion} we get 
                                                    $$
                                                    \alpha_{j+1} = 2 - \dfrac{1}{\alpha_{j+2}} = 1 + \dfrac{1}{k-j},
                                                    $$ 
                                                     as wanted. 
                                                       \qed 
                                                       
\pf{of Theorem \ref{thm-main}} The result follows from Lemmas \ref{lemma-xk-arb-large} through \ref{lemma-alphaj-smallest}.
Precisely, by Lemma \ref{lemma-xk-arb-large} variable $x_k$ can be arbitrarily large in a strictly feasible solution 
of \eqref{problem-sdpprime}. 
By Lemma
\ref{lemma-derive-poly}  we derive the polynomial inequalities  \eqref{eqn-pj}. 
From these in Lemma \ref{lemma-xj>=dj+1} we derive the clean inequalities 
\eqref{eqn-xjgeqdj+1} via the recursion \eqref{eqn-alpha-j-recursion}. 

From the  recursion \eqref{eqn-alpha-j-recursion} it directly follows that all $\alpha_{j+1}$ are at most $2.$ The lower bound 
on the $\alpha_{j+1}$ is proved as follows: by Lemma \ref{lemma-alphaj-monotone} the $\alpha_{j+1}$ are monotone functions of the tail-indices $t_{j+1}.$ On the other hand,   $t_{j+1} \geq j+2$ for all $j$ by Lemma \ref{lemma-tailindex}
and when $t_{j+1} = j+2 \,$ for all $j, \, $ then  by Lemma \ref{lemma-alphaj-smallest} we have $\alpha_{j+1} = 1 + 1/(k-j).$ The proof is now complete. \qed

\subsection{Computing the exponents by Fourier-Motzkin elimination} 
\label{subsection-fourier-motzkin}

The recursion  \eqref{eqn-alpha-j-recursion}  gives a convenient way to compute the 
$\alpha_j$ exponents. Equivalently, we can compute the $\alpha_j$ via the 
 well known Fourier-Motzkin elimination algorithm, designed for linear inequalities; this is an interesting 
 contrast, since SDPs are  highly nonlinear. 

We do this as follows. If polynomial $p_j$ is of type 1, then we suppress the lower order terms to get 
\begin{equation} 
x_j x_{t_{j+1}} \geq \const x_{j+1}^2,
\end{equation}
see the last inequality in \eqref{eqn-jerry-3}. If polynomial $p_j$ is of type 2,  then we similarly
 suppress the lower order terms to deduce 
 		\begin{equation}
 		x_j \geq \const  x_{j+1}^2. 
 		\end{equation}
 After this, using that $x_1, \dots, x_k$ are all positive, we  rewrite the inequalities in terms 
 of $y_j := \log_2 x_j$ for all $j, \, $ then eliminate variables. For example, from the inequalities 
 		\eqref{eqn-mild-quadratic} we deduce 
 		\begin{equation} \label{eqn-mild-quadratic-log} 
 		\begin{array}{rcl} 
 		y_1 + y_3   & \geq & 2 y_2   \\ 
 		y_2 + y_4   & \geq&  2 y_3 \,\,\, \\ 
 		y_3  & \geq &  2 y_4.   \,\,\, 
 		\end{array} 
 		\end{equation}
 		We add $\frac{1}{2}$ times the last inequality in \eqref{eqn-mild-quadratic-log} to 
 		the middle one  to get 
 		\begin{equation} \label{eqn-new1} 
 		y_2  \geq \frac{3}{2} y_3.
 		\end{equation}
 		We then add $\frac{2}{3}$ times \eqref{eqn-new1} to the first inequality in  \eqref{eqn-mild-quadratic-log} 	to get 
 		\begin{equation}   \label{eqn-new2} 
 		y_1 \geq \frac{4}{3} y_2.
 		\end{equation} 
 	 Finally,  \eqref{eqn-new1}, \eqref{eqn-new2} and the last inequality in \eqref{eqn-mild-quadratic-log}   translate back to the   inequalities 
 		 		\eqref{eqn-mild-quadrati-corollary}.

 		\section{When we do not even need a change of variables}
 		\label{section-polyopt} 
 		
 		As we previously discussed, the linear change of variables $x \leftarrow Mx$ 
 		 is necessary to obtain a Khachiyan type hierarchy among the variables. 
 		 		 		Nevertheless, in this section we show a natural SDP in which 
 		 		 		a Khachiyan type hierarchy occurs even without a change of variables; more precisely, the SDP is in the form of \eqref{problem-sdpprime}. For completeness, we also revisit O' Donnell's example from   \cite{o2017sos}, and show 
 		 	that the SDP therein is also in the regular form of \eqref{problem-sdpprime}. 
 	
Given a univariate polynomial of even degree $f(x) = \sum_{i=0}^{2n} a_i x^i$ with $a_{2n} > 0,$  
we 
consider  the problem of minimizing $f$ over $\rad{}.$ We write this problem as
	\begin{equation} \label{problem-polyopt} 
	\begin{array}{rlcl} 
	\sup & \lambda \\
	s.t.   & f  - \lambda & \geq & 0.
	\end{array}
	\end{equation}
	We will show that in the natural SDP formulation of \eqref{problem-polyopt} exponentially large variables appear naturally, although here by ``exponentially large" we only mean in magnitude, not in size.
	
		Since $f - \lambda$ is also a univariate polynomial, it is nonnegative if and only if   it is a sum of squares (SOS), that is, iff 
		$f - \lambda = \sum_{i=1}^t g_i^2$ for a positive integer $t$ and polynomials $g_i$ \footnote{However, there are multivariate polynomials that are nonnegative, but {\em not} SOS.}. 
			Define the vector of monomials 
			$$
			z = (1, x, x^2, \dots, x^n)^\top.
			$$
			Then $f-\lambda$ is SOS if and only if (see \cite{lasserre2001global, nesterov2000squared, parrilo2003semidefinite, shor1987class}) 
		$
		f - \lambda = zz^\top \bullet Q
		$
		for some $Q \succeq 0.$ 
		We then match monomials in $f - \lambda$ and  $zz^\top \bullet Q$ and translate \eqref{problem-polyopt} into the SDP
	
	\begin{equation} \label{problem-polyopt-sdp} 
	\begin{array}{rrcl} 
	\max & - A_0 \bullet Q \\
	s.t.     & A_i \bullet Q & = & a_i \,\, \text{for} \, i=1, \dots, 2 n \\
	           &                   Q & \in & \psd{n+1}.
	\end{array} 
	\end{equation}
	Here for all $i \in \{0, 1, \dots, 2 n \}$ the $(k, \ell)$ element of the matrix $A_i$ is $1$ if $k + \ell = i+2$ for some $k, \ell \in \{1, \dots, n+1\}, \,$ 
	and all other entries of $A_i$ are zero.
	
For positive integers $k$ and $\ell, \,$ let us define $E_{k \ell}$ as the $(k, \ell)$th unit matrix, whose $(k, \ell)$ and $(\ell,k)$ entries are $1, \, $ and the rest are zero.

	\begin{Lemma} \label{lemma-polyopt} 
		After permuting and renaming variables, the constraints of the dual problem of \eqref{problem-polyopt-sdp} can be written  as 
			\begin{equation} \label{eqn-polyopt-dual-constraints-2} 
			\sum_{i=1}^{2n} x_i A_i' + E_{n+1,n+1} \succeq 0,
			\end{equation}
			with $
			A_i' = \sum_{k + \ell = 2i} E_{k \ell} \; \text{for} \; i=1, \dots, n.
			$
			Thus $(A_1', A_2', \dots, A_n')$ is a regular facial reduction sequence, and the constraint set \eqref{eqn-polyopt-dual-constraints-2} is in the form 
			of \eqref{problem-sdpprime}, with $k=n.$ 
			\qed
	\end{Lemma}

 Before we prove Lemma \ref{lemma-polyopt}, we illustrate it. For that, suppose $n=3.$ Then by Lemma \ref{lemma-polyopt} 
		the constraints \eqref{eqn-polyopt-dual-constraints-2} look like 
		$$
		x_1	\underbrace{\begin{pmatrix}
			1       & \pha{0}   & \pha{0}  & \pha{0}    \\
			\pha{0}    & {0}   & \pha{0}  & \pha{0}    \\
			\pha{0}     & \pha{0}   & {0}  & \pha{0}    \\
			\pha{0}     & \pha{0}   & \pha{0}  & {0} 
			\end{pmatrix}}_{A_1'}  + x_2		\underbrace{\begin{pmatrix}
			0       & \pha{0}   & 1  & \pha{0}    \\
			\pha{0}    & {1}   & \pha{0}  & \pha{0}    \\
			1      & \pha{0}   & {0}  & \pha{0}    \\
			\pha{0}     & \pha{0}   & \pha{0}  & {0} 
			\end{pmatrix}}_{A_2'}  + x_3	\underbrace{\begin{pmatrix}
			0       & \pha{0}   & \pha{0}  & \pha{0}    \\
			\pha{0}    & {0}   & \pha{0}  & 1    \\
			\pha{0}     & \pha{0}   & {1}  & \pha{0}    \\
			\pha{0}     & 1   & \pha{0}  & {0} 
			\end{pmatrix}}_{A_3'}  + \sum_{i=4}^6 x_i A_i'  + E_{4,4} \succeq 0.
		$$
		Note that if we delete the term $\sum_{i=4}^6 x_i A_i'$ from this system, we obtain a smaller version of 
		our previously discussed 
		problem 	\eqref{problem-mild-sdp} 	 (with three variables rather than four).
	
	\pf{of Lemma \ref{lemma-polyopt}:} 
	 The dual problem of \eqref{problem-polyopt-sdp} is 
	\begin{equation} \label{problem-polyopt-sdp-dual} 
	\begin{array}{rrcl} 
	\min  & \sum_{i=1}^{2n} a_i y_i \\ 
	s.t.     & \sum_{i=1}^{2n} y_i  A_i + A_0 & \succeq & 0,
	\end{array} 
	\end{equation}
	whose constraints  can be written as 
	$$
	\begin{pmatrix}
	1      & y_1 & y_2 & \dots &   y_n \\
	y_1       & y_2 &        & \dots & y_{n+1} \\
	y_2       &        &        &  \dots    & y_{n+2}       \\
	\vdots &         &        & \ddots  &   \vdots                     \\
	y_n      & y_{n+1}        &   y_{n+2}     &   \dots   & y_{2n}              
	\end{pmatrix} \succeq 0. 
	$$
	Permuting rows and columns, this is equivalent to
	\begin{equation} \label{eqn-polyopt-dual-constraints} 
	\begin{pmatrix}
	y_{2n}          & y_{2n-1}  & y_{2n-2}  & \dots &   y_n \\
	y_{2n-1}       & y_{2n-2} &        & \dots & y_{n-1} \\
	y_{2n-2}       &        &        &  \dots    & y_{n-2}       \\
	\vdots &         &        & \ddots  &   \vdots                     \\
	y_n      & y_{n-1}        &   y_{n-2}     &   \dots   & 1           
	\end{pmatrix} \succeq 0.
	\end{equation} 
	Let us rename the variables so the even numbered ones come first, and the rest come afterwards,
	as 
	\begin{equation} \label{eqn-rename-y-x} 
	\begin{array}{rclrclcrcl} 
	x_1 \! & \! := \! & \! y_{2n}, \, & x_2 \! & \! := \! & \! y_{2n-2}, & \dots &  x_n \! & \! := \! & \! y_2; \\
	x_{n+1} \! & \! := \! & \! y_{2n-1}, & x_{n+2} \! & \! := \! & y_{2n-3}, & \dots &  x_{2n} \! & \! := \! & \! y_1. 
	\end{array} 
	\end{equation} 
	Then the constraints \eqref{eqn-polyopt-dual-constraints} become as required in \eqref{eqn-polyopt-dual-constraints-2} with $(A_1', A_2', \dots, A_n')$ being a regular facial reduction sequence. Finally, an argument just like the one after  
	 Example \ref{example-khanch-sdp} shows that the singularity degree of 
	 $\{ \, Y \succeq 0 : \, A_i' \bullet Y = 0 \; \text{for} \; i=1, \dots, 2 n \, \}$ is $n.$ For this last part, we leave the details to the reader. 
	\qed

	We next claim that for all feasible solutions of \eqref{eqn-polyopt-dual-constraints-2} the inequalities 
		\begin{equation} \label{eqn-xjxj+2} 
	x_j x_{j+2} \geq x_{j+1}^2 \, \text{for} \, j=1, \dots, n-2; \; {\rm and} \; x_{n-1} \geq x_n^2
	\end{equation} 
	hold. Indeed, we can derive these by following the proof  of Lemma \ref{lemma-derive-poly}, since the tail-indices (cf. Definition  \ref{definition-tailindex}) are  
	$t_{j+1} = j+2$ for $j=1, \dots, n-1.$ Further, now the $``\delta"$ terms that appear in Lemma \ref{lemma-derive-poly} are all zero, so we do not have to worry about strict feasibility, nor about ``making $x_k$ large." 	Hence by Lemma \ref{lemma-alphaj-smallest} we deduce that 
	\begin{equation} \label{kramer} 
	x_j  \geq x_{j+1}^{\alpha_{j+1}} \, \text{for} \, j=1, \dots, n-1
	\end{equation} 
	hold, where $\alpha_{j+1} = 1 + 1/(n-j)$ for all $j.$ From these inequalities we then derive \footnote{Precisely, by straightforward induction we get  $x_1 \geq x_{j+1}^{n/(n-j)}$ for $j=1, \dots, n-1.$ } 
	\begin{equation} \label{george} 
		x_1 \geq x_n^n.
	\end{equation}

	We invite the reader to try a simpler alternative derivation of the inequalities above: 
	first,  from the order two subdeterminants in 
	\eqref{eqn-polyopt-dual-constraints-2} directly deduce the inequalities \eqref{eqn-xjxj+2}; and second, from 	\eqref{eqn-xjxj+2}
	eliminate variables to get the inequalities in  \eqref{kramer}. 
	
	We translate the  inequalities in \eqref{kramer} and \eqref{george} back to the original $y_j$ variables (using  the correspondence \eqref{eqn-rename-y-x}), and obtain the following result: 
\co{	We translate the  inequalities in the $x$ space to the space of the $y$ variables, 
	and obtain 
	the following result:} 
	\begin{Theorem} \label{thm-polyopt} 
		Suppose that $y \in \rad{2n}$ is feasible in \eqref{problem-polyopt-sdp-dual}. Then the following hold:
	\begin{eqnarray} \label{thm-polyopt-1} 
			y_{2(n-j+1)} &  \geq &	y_{2(n-j)}^{1 + 1/(n-j)}  \,  \text{for} \, j=1, \dots, n-1 \\ \label{thm-polyopt-2} 
				y_{2n} & \geq & y_2^n.
	\end{eqnarray}
\qed
	\end{Theorem}
	We next connect Theorem \ref{thm-polyopt} to other results in the literature.
	
	First, Theorem \ref{thm-polyopt} complements a result of Lasserre \cite[Theorem 3.2]{lasserre2001global}, which states the following:
	if $\bar{x}$ minimizes the  polynomial $f(x)$ then
	$$
	(y_1,y_2,\dots,y_{2n})=(\bar{x},\bar{x}^2,\dots,\bar{x}^{2n})
	$$
	is optimal in \eqref{problem-polyopt-sdp-dual}. On the one hand, Theorem \ref{thm-polyopt} states bounds on  all feasible solutions,  on the other hand, it does not specify an optimal solution. 
	
Second, the matrix in formula \eqref{eqn-polyopt-dual-constraints} is a Hankel matrix, i.e., 
its elements along the reverse diagonals are constant
 \footnote{The matrix in \eqref{eqn-polyopt-dual-constraints} has a ``1" in the lower right corner, but this can always be achieved by a  normalization.}. 
  Theorem \ref{thm-polyopt} compares the even numbered elements of this matrix,
no matter what the odd numbered elements 
$y_{2n-1}, y_{2n-3}, \dots, y_3$ are. Thus, it is related to  recent results  of Choi and Jafari  \cite{choi2016positive} on partially defined matrices,
which can be completed to be Hankel, and positive (semi)definite. 
 \co{Third, Theorem \ref{thm-polyopt} states that in any psd Hankel matrix the elements must be very different. Thus, in a broader sense
 they complement  results of Beckermannn and Townsend \cite{beckermann2017singular} and Tyrtishnikov \cite{tyrtyshnikov1994bad}  
about the exponentially ill conditioned nature of psd Hankel matrices.
}
	
		For completeness, we next revisit an example of O' Donnell in \cite{o2017sos}, and show how the SDP that appears 
		in there is in the regular form of \eqref{problem-sdpprime}.

		\begin{Example}
		We are given the 
	polynomial with $2n$ variables 
	$$
	p(x,y) = p(x_1, \dots, x_n, y_1, \dots, y_n) = x_1 + \dots + x_n - 2 y_1
	$$
	and the set $K$ defined by the equations 
		$$
	\begin{array}{rclcrclcrclcrclrcl}
	2 x_1 y_1 \!\! & \! = \! & \!\! y_1, & & 2x_2 y_2 \!\! & \! = \! & \!\! y_2, &  \dots & 2 x_{n-1}  y_{n-1}  \!\! & \! = \! & \!\! y_{n-1},    & 2 x_n y_n \!\! & \! = \! & \!\! y_n, \\
	x_1^2 \!\! & \! = \! & \!\! x_1,      & & x_2^2 \!\! & \! = \! & \!\! x_2,       & \dots & x_{n-1}^2 \!\! & \! = \! & \!\! x_{n-1}, &  x_n^2 \!\! & \! = \! & \!\! x_n, \\ 
	y_1^2 \!\! & \! = \! & y_2,            & & y_2^2 \!\! & \! = \! & y_3,            &  \dots & y_{n-1}^2 \!\! & \! = \! & y_{n},  &  y_n^2 \!\! & \! = \! & 0. 
	\end{array}
	$$
Note that in the description of $K$ the very last constraint $y_n^2 = 0$ breaks the pattern seen in the previous $n-1$ columns. 	We ask the following question:
	\begin{itemize}
		\item Is $p(x,y) \geq 0$ for all $(x,y) \in K$?
	\end{itemize}
	The answer is clearly yes, since for all $(x,y) \in  K$ we have $x_1, \dots, x_n  \in \{0,1\}$ and 
	$y_1 = \dots = y_n = 0.$ 
	
	On the other hand, the sum of squares  procedure verifies the ``yes" answer as follows.
	Let 
	$$
	z = (1, x_1, \dots, x_n, y_1, \dots, y_n)^\top
	$$
	be a vector of monomials. To certify that $p(x,y)$ is nonnegative over $K, $ we seek $\lambda, \mu, \nu \in \rad{n}$ and $Q \succeq 0$ such that 
	\begin{equation} \label{eqn-p-Q-lambda} 
	\begin{array}{rcl}
	p(x,y) =z^\top Qz & + & \lambda_1(2x_1y_1-y_1)+\mu_1(x_1^2-x_1)+\nu_1(y_1^2-y_2) \\
	& + & \lambda_2(2x_2y_2-y_2)+\mu_2(x_2^2-x_2)+\nu_2(y_2^2-y_3) \\
	&\vdots& \\
	& + & \lambda_n(2x_ny_n-y_n)+\mu_n(x_n^2-x_n)+\nu_n(y_n^2-0).
	\end{array}
		\end{equation}
	Indeed, if we succeed and find such $\lambda, \mu, \nu, Q, \, $  then $p(x,y) = z^\top Q z \geq 0$  for all $(x,y) \in K.$
		
	 The polynomials (in the $x_i$ and $y_i$)  on the two sides of 
		\eqref{eqn-p-Q-lambda}  are equal exactly when all their coefficients are equal. 
		Thus, matching coefficients in 	\eqref{eqn-p-Q-lambda}  on the left and right hand sides, \mbox{O' Donnell}  \cite{o2017sos} showed that any 
	$Q$ feasible  in \eqref{eqn-p-Q-lambda} looks like 
	\begin{equation} \label{eqn-Q-lookslikethis} 
	Q \, = \, \left(
	\begin{array}{ccccc|ccccc|c}
	u_1 & 0 & \dots & 0 & 0 & -u_2 & \dots & 0 & 0 & 0 & 0 \\
	0 & u_2 & \dots & 0 & 0 & \vdots & \ddots & \vdots & \vdots & \vdots & \vdots \\
	\vdots & \vdots & \ddots & \vdots & \vdots & 0 & \dots & -u_{n-1} & 0 & 0 & 0 \\
	0 & 0 & \dots & u_{n-1} & 0 & 0 & \dots & 0 & -u_n & 0 & 0 \\ 
	0 & 0 & \dots & 0 & u_{n} & 0 & \dots & 0 & 0 & -2 & 0 \\ \hline
	-u_2 & \dots & 0 & 0 & 0 & 1 & \dots & 0 & 0 & 0 & 0 \\
	\vdots & \ddots & \vdots & \vdots & \vdots & \vdots & \ddots & \vdots & \vdots & \vdots & \vdots \\
	0 & \dots & -u_{n-1} & 0 & 0 & 0 & \dots & 1 & 0 & 0 & 0 \\
	0 & \dots & 0 & -u_n & 0 & 0 & \dots & 0 & 1 & 0 & 0 \\
	0 & \dots & 0 & 0 & -2 & 0 & \dots & 0 & 0 & 1 & 0 \\ \hline
	0 & \dots & 0 & 0 & 0 & 0 & \dots & 0 & 0 & 0 & 0
	\end{array}
	\right)
	\end{equation}
	for suitable $u_1, \dots, u_n.$ 
	Looking at $2 \times 2$ subdeterminants of $Q$ we see that 
	 the $u_i$ satisfy 
	\begin{equation} \label{eqn-odonnell} 
	u_1 \geq u_2^2, u_2 \geq u_3^2, \dots, u_{n-1} \geq u_n^2, u_n \geq 4. 
	\end{equation} 
which is the same as \eqref{problem-khachiyan}, except we replaced the constant $2$ by $4.$

	Recall that $E_{k \ell}$ is the unit matrix   
	in which the $(k, \ell)$ and $(\ell, k)$ entries are $1$ and the rest zero. Define 
	$$
	A_1' = E_{11}, A_i' = E_{ii} - E_{i-1, n+i-1} \; \text{for} \; i =2, \dots, n.
	$$
	Then any $Q$ feasible in \eqref{eqn-p-Q-lambda} is written as 
	\begin{equation}  \label{eqn-Q} 
	Q = u_1 A_1' + u_2 A_2' + \dots + u_n A_n' + B' \succeq 0, 
	\end{equation} 
	for a suitable $B'$ (precisely, $B'     =  - 2 E_{n, 2n} + \sum_{i=n+1}^{2n} E_{ii}$). 
	
	We see that $(A_1', \dots, A_n')$ is a regular facial reduction sequence, 
	thus the system \eqref{eqn-Q} is in 
	the regular form of \eqref{problem-sdpprime}.

		\end{Example}

	Note that \eqref{eqn-Q-lookslikethis} arises by concatenating $2 \times 2$ psd blocks of the form 
	\eqref{eqn-khach-2by2}, then permuting rows and columns. In other words, 
	\eqref{eqn-Q-lookslikethis} is the exact representation of \eqref{problem-khachiyan} (apart from the constant $2$ being replaced by $4$ and the $u_i$ being negated in the offdiagonal positions), that we discussed after Example \ref{example-khanch-sdp}. 
	
Among followup papers of O' Donnell  \cite{o2017sos} we should mention the work of Raghavendra and Weitz \cite{raghavendra2017bit} 
which gave SDPs which also have  a sum-of-squares origin, and 
exponentially large size solutions. It would be interesting to see whether those SDPs are also in the regular  form of 
\eqref{problem-sdpprime}.

\section{Conclusion}
\label{section-conclusion} 

Khachiyan's SDP is a classical pathological problem in which the size of any feasible solution is exponential in the number of variables.
Here we showed that Khachiyan's SDP is far from being an isolated example:  in any strictly feasible SDP 
a linear transformation induces a Khachiyan type hierarchy among a subset of the variables, 
and large size solutions. 
The number of variables in the hierarchy and ``how large" they get, depends on the singularity degree of a dual problem. 
Further, such a hierarchy and large solutions naturally appear in SDPs that come from sum-of-squares optimization,  
without any change of variables. 

We also studied how to represent large solutions of SDPs in polynomial space. Our main tool was
the regularized semidefinite program \eqref{problem-sdpprime}. 
If \eqref{problem-sdp} and \eqref{problem-sdpprime} are strictly feasible, then 
in the latter 
we can verify that a strictly feasible solution exists, without computing the actual values of the ``large" variables $x_1, \dots, x_k:$ see Figure \ref{figure-add-xk}.  Further, SDPs that arise from polynomial optimization (Section \ref{section-polyopt}) and the SDP that represents \eqref{problem-khachiyan} are naturally in the 
form of \eqref{problem-sdpprime}. Hence in these SDPs we can also certify large solutions without computing their actual values. 

Several questions remain open. For example, what can we say about large solutions 
in semidefinite programs that are not strictly feasible? The discussion after Example \ref{example-khanch-sdp} shows that we do not have a complete answer yet. 

Also, recall that we transform \eqref{problem-sdp} into \eqref{problem-sdpprime} by a linear change of variables
(equivalent to operations \eqref{exch} and \eqref{trans} in Definition \ref{definition-reform}) and a similarity transformation (operation \eqref{rotate} in Definition \ref{definition-reform}). The latter has no effect on how large the variables are.
We are thus led to the following question: are {\em all}  \ SDPs with exponentially large solutions  in the form of \eqref{problem-sdpprime} (perhaps after a similarity transformation)?  In other words, can we {\em always} certify large size solutions in SDPs using a regular facial reduction sequence? Answering this question would help us answer the greater question: can we decide feasibility of SDPs in polynomial time?

\co{
Some natural open questions remain. For example, we did not address the issue of how large 
the entries in the matrices in \eqref{problem-sdpprime} can become. Another question that may help us solve deciding SDP feasibility is whether {\em any} SDP with exponentially large solutions is in the form of \eqref{problem-sdpprime}. In other words, can we certify large size solutions in every SDP using a regular facial reduction sequence?
}

{\bf Acknowledgement} We are very grateful for helpful discussions to Richard Rim\'anyi, in particular, to suggestions  that led to the proof of Lemma \ref{lemma-alphaj-monotone}. 
We are also very grateful to 
the anonymous referees and  to Andrew B. Nobel whose suggestions helped us to improve the paper. 
We thank Alex Townsend for pointing out reference \cite{choi2016positive}.
We gratefully acknowledge the support of National Science Foundation, award DMS-1817272.

 \bibliographystyle{plain}
 \bibliography{mysdpMelody}

\end{document}